\documentclass[aps,prd,10pt,notitlepage,nofootinbib,superscriptaddress,showkeys,showpacs]{revtex4-1}

\usepackage{amsfonts, amsmath, amssymb, amsthm}
\usepackage{bm}
\usepackage{color}
\usepackage{stmaryrd}
\usepackage{graphicx}
\usepackage[babel=true]{csquotes}
\usepackage{amsmath}
\usepackage{bbm}
\usepackage{epstopdf}
\DeclareMathOperator{\Tr}{Tr}

\newcommand{\be}{\begin{equation}}
\newcommand{\ee}{\end{equation}}
\newcommand{\bea}{\begin{eqnarray}}
\newcommand{\eea}{\end{eqnarray}}

\newcommand{\cM}{\mathcal{M}}

\newcommand{\cB}{\mathcal{B}}
\newcommand{\cF}{\mathcal{F}}

\newcommand{\cE}{\mathcal{E}}
\newcommand{\cV}{\mathcal{V}}
\newcommand{\cW}{\mathcal{W}}
\newcommand{\cG}{\mathcal{G}}

\newcommand{\cT}{\mathcal{T}}
\newcommand{\cP}{\mathcal{P}}

\newcommand{\cI}{\mathcal{I}}
\newcommand{\cJ}{\mathcal{J}}

\newcommand{\bM}{\mathbb{M}}
\newcommand{\bG}{\mathbb{G}}
\newcommand{\bB}{\mathbb{B}}
\newcommand{\bC}{\mathbb{C}}
\newcommand{\bN}{\mathbb{N}}
\newcommand{\bW}{\mathbb{W}}

\newcommand{\un}{\mathbbm{1}}

\newcommand{\fsig}{\bm{\sigma}}

\newcommand{\hI}{\widehat{\cI}}

\newcommand{\prf}{{\noindent \bf Proof\; \; }}

\DeclareMathOperator{\col}{col}
\DeclareMathOperator{\opt}{opt}

\newtheorem{lemma}{Lemma}
\newtheorem{definition}{Definition}
\newtheorem{theorem}{Theorem}
\newtheorem{coroll}{Corollary}
\newtheorem{prop}{Proposition}

\begin{document}

\title{\Large Colored triangulations of arbitrary dimensions are stuffed Walsh maps}

\author{{\bf Valentin Bonzom}}\email{bonzom@lipn.univ-paris13.fr}
\affiliation{LIPN, UMR CNRS 7030, Institut Galil\'ee, Universit\'e Paris 13, Sorbonne Paris Cit\'e, 99, avenue Jean-Baptiste Cl\'ement, 93430 Villetaneuse, France, EU}

\author{{\bf L.~Lionni}}\email{luca.lionni@th.u-psud.fr}
\affiliation{Laboratoire de Physique Th\'eorique, CNRS UMR 8627, Universit\'e Paris XI, 91405 Orsay Cedex, France, EU}
\affiliation{LIPN, UMR CNRS 7030, Institut Galil\'ee, Universit\'e Paris 13, Sorbonne Paris Cit\'e, 99, avenue Jean-Baptiste Cl\'ement, 93430 Villetaneuse, France, EU}

\author{{\bf V.~Rivasseau}}\email{rivass@th.u-psud.fr}
\affiliation{Laboratoire de Physique Th\'eorique, CNRS UMR 8627, Universit\'e Paris XI, 91405 Orsay Cedex, France, EU}
\affiliation{Perimeter Institute for Theoretical Physics, 31 Caroline St. N, ON N2L 2Y5, Waterloo, Canada}

\begin{abstract}
Regular edge-colored graphs encode colored triangulations of pseudo-manifolds. Here we study families of edge-colored graphs built from a finite but arbitrary set of building blocks, which extend the notion of $p$-angulations to arbitrary dimensions. We prove the existence of a bijection between any such family and some colored combinatorial maps which we call stuffed Walsh maps. Those maps generalize Walsh's representation of hypermaps as bipartite maps, by replacing the vertices which correspond to hyperedges with non-properly-edge-colored maps. This shows the equivalence of tensor models with multi-trace, multi-matrix models by extending the intermediate field method perturbatively to any model. We further use the bijection to study the graphs which maximize the number of faces at fixed number of vertices and provide examples where the corresponding stuffed Walsh maps can be completely characterized.
\end{abstract}


\maketitle

\section*{Introduction}

Regular, bipartite, edge-colored graphs are connected, bipartite graphs with vertices of degree $(D+1)$ such that the $(D+1)$ incident edges on each vertex all have a distinct color from the set $\{0, \dotsc, D\}$ (see an example in Figure \ref{fig:VacFeynm}). Throughout this article, we refer to them simply as edge-colored graphs. By duality, they represent colored triangulations of pseudo-manifolds in dimension $D$ \cite{ColoredReview}. Each vertex represents a $D$-simplex whose boundary simplices of dimension $D-1$ are colored since they are dual to the edges of the graph. Two $D$-simplices can be glued by identifying two of their boundary $(D-1)$-simplices if they have the same color. To specify the identification of those $(D-1)$-simplices, one looks at the $(D-2)$-simplices which are colored by a pair of colors and identify two of them if they have the same colors. This reasoning is repeated down to the vertices of the simplices. This shows that an edge-colored graph determines a unique gluing of $D$-simplices.

The vertices of a colored triangulation are dual to the edge-colored subgraphs regular of degree $D$ carrying all the colors but one. We will be interested in those which have the colors $1, \dotsc, D$ and call them \emph{bubbles}. A bubble is a regular, bipartite, edge-colored graph with $D$ colors, dual to a triangulation of a pseudo-manifold in dimension $D-1$ and gives a $D$-dimensional object with boundary using topological cones \cite{Uncolored}. In two dimensions, bubbles are simply bipartite cycles with edges alternating the colors 1 and 2. Taking the topological cone turns them into quadrangles, hexangles, etc. We can thus think of a bubble as a building block for colored triangulations of pseudo-manifolds in arbitrary dimensions, generalizing the triangles, quadrangles, etc. which are the two-dimensional building blocks and think of edge-colored graphs with fixed bubbles as a higher-dimensional extension of combinatorial maps with prescribed face degrees. We denote $\bG_0(\{\cB_\alpha\}_{\alpha\in A})$ the set of edge-colored graphs with $D+1$ colors whose bubbles are from the finite set of bubbles $\{\cB_\alpha\}_{\alpha\in A}$.

A special attention has been recently paid to colored graphs with quartic bubbles (four vertices), first with the so-called quartic melonic bubbles \cite{BeyondPert} (see Figure \ref{fig:BubbleExamples}) and then with generic quartic bubbles \cite{ArbitraryQuartic}. The reason\footnote{The reason is in fact stronger than the bijection we describe. It is due to a non-perturbative equivalence which holds between tensor integrals and matrix integrals, which are respectively perturbative generating functions of colored graphs and colored maps.} is that in those cases, there is an extension of Tutte's bijection between bipartite quadrangulations and generic maps which applies. It gives a bijection between edge-colored graphs with quartic bubbles and generic maps whose edges can be labeled by a color or a subset of colors.

This bijection has been really fruitful. In particular, it was used in \cite{DSDartois} to perform a double scaling limit in the quartic melonic case (this is a scaling limit of the generating function which picks up the most singular contributions for all values of the difference between the number of faces and the number of vertices).

For non-quartic bubbles, very few results exist, see for instance \cite{Meanders, UnitaryInt} where graphs with a single bubble, and for specific bubbles, are studied, focusing on those which maximize the number of faces. In a few instances, a non-direct approach was developed. It was possible to rely on results derived for quartic bubbles, typically via the bijection, and extend them to more generic (but not arbitrary) bubbles. In particular, this strategy was applied in \cite{DSSD} to derive the double scaling limit for generic melonic bubbles and in \cite{Melonoplanar} for the analysis of \enquote{melono-planar} graphs. Extending results from the quartic case to more generic bubbles $\cB$ was possible via a universality argument. The bubbles $\cB$ can be represented as gluings of quartic bubbles so that the corresponding set of colored graphs $\bG_0(\cB)$ forms a subfamily of those built from quartic bubbles. If $\cB$ is \enquote{close enough} to the quartic ones, it may lead to the same combinatorial universality class.

However, this argument cannot pass to fully generic bubbles. If $\cB$ becomes complicated, then $\bG_0(\cB)$ will span a very special subset of the graphs with quartic bubbles which is difficult to identify and study. Therefore the universality argument becomes out of reach. This problem suggests that instead of the non-direct approach which goes through the quartic model for which a bijection with maps exists, one should look for a direct approach and find a bijection on $\bG_0(\cB)$. That bijection is the main result of this article. Yet, the fact that all bubbles can be represented as gluings of quartic melonic bubbles (a statement proved for the first time in the present article) is the key to extend the bijection from quartic bubbles to arbitrary ones.

In the quartic case, the bijection works as follows. Quartic bubbles can be represented as colored edges, where each endpoint of an edge represents a pair of vertices of the quartic bubble. Since gluing bubbles via edges of color 0 forms cycles, the latter are then represented as vertices locally embedded in the plane (this will be detailed in Section \ref{sec:Review}), thereby leading to (non-properly-)edge-colored maps. A generic bubble however has more than two pairs of vertices, hence through the bijection the \enquote{edges} should have more than two ends. This is reminiscent of hyperedges of hypermaps. Following Walsh \cite{Walsh}, a hypermap can be represented as a bipartite map, sometimes known in this context as a Walsh map. Indeed, a hyperedge with $n$ ends (called bits in \cite{Walsh}) can be represented as a star centered on a white vertex incident to $n$ half-edges. Those half-edges can then be glued to half-edges incident to ordinary, say black, vertices and produce bipartite maps.

However, a bubble, with its colored edges, has in general more structure than a hyperedge. We will show that this structure can always be represented through edge-colored maps. Then, instead of Walsh's white vertices representing hyperedges, we will insert those edge-colored maps. This process is a stuffing of Walsh maps with specific edge-colored maps. It is thus a bipartite, colored extension of the stuffed maps introduced in \cite{Stuffed}, which we will call \emph{stuffed Walsh maps}.

After setting up definitions and notations and reviewing the bijection in the quartic melonic case in Section \ref{sec:Review}, we present the generalized bijection between $\bG_0(\cB)$ and stuffed Walsh maps in Section \ref{sec:Bijection}. Although we focus for convenience on the case of a single type of bubble $\cB$, it should be clear that the bijection generalizes to $\bG_0(\{\cB_\alpha\}_{\alpha\in A})$ for a finite set of bubbles $\{\cB_\alpha\}_{\alpha\in A}$. We then use the bijection to study the \emph{dominant edge-colored graphs}, i.e. those which maximize the number of faces at fixed number of vertices. Section \ref{sec:DominantMaps} introduces a new approach to analyze this issue, based on the genera and circuit ranks of monochromatic submaps of stuffed Walsh maps. While characterizing the dominant maps remains an open issue in general, we have provided examples in Section \ref{sec:Applications} for which the identification of dominant maps can be fully carried out. This includes bubbles for which the problem of the dominant graphs had not been solved before.

Finally, in Section \ref{sec:MatrixModel}, we present a multi-matrix integral which generates stuffed Walsh maps. This is expected since multi-trace matrix models generate stuffed maps as shown in \cite{Stuffed}. As the generating function of regular edge-colored graphs takes the form of an integral over a tensor, known as a tensor model \cite{Uncolored}, our results mean that tensor models are (perturbatively) equivalent to some multi-trace multi-matrix models. This is proved directly in Section \ref{sec:MatrixModel} using matrix and tensor integrals.

As the bijection relies on representing bubbles as edge-colored maps, we provide details on the admissible maps in Appendix \ref{sec:EdgeColoredMaps}, using the same techniques as in Section \ref{sec:DominantMaps}. In Appendix \ref{HubStrat},  we reduce the matrix model of Section \ref{sec:MatrixModel} in the special case of the quartic melonic bubbles to the matrix model of \cite{BeyondPert}.

\section{Edge-colored graphs with quartic melonic bubbles and edge-colored maps} \label{sec:Review}

A \emph{closed $\Delta$-edge-colored} graph is a connected, bipartite, regular graph in which each edge carries a color from $\{1, \dotsc, \Delta\}$ and such that each vertex is incident to exactly $\Delta$ edges of distinct colors. An \emph{open} $\Delta$-colored graph may have some vertices with less than $\Delta$ incident edges (all with distinct colors).

Let $D\geq 2$ be an integer. A \emph{bubble} is a closed $D$-colored graph with colors in $[D]\equiv \{1, \dotsc, D\}$ and we denote $\bB_p$ the set of bubbles with $2p$ vertices ($p$ white and $p$ black vertices).

\begin{figure}
\includegraphics[scale=0.4]{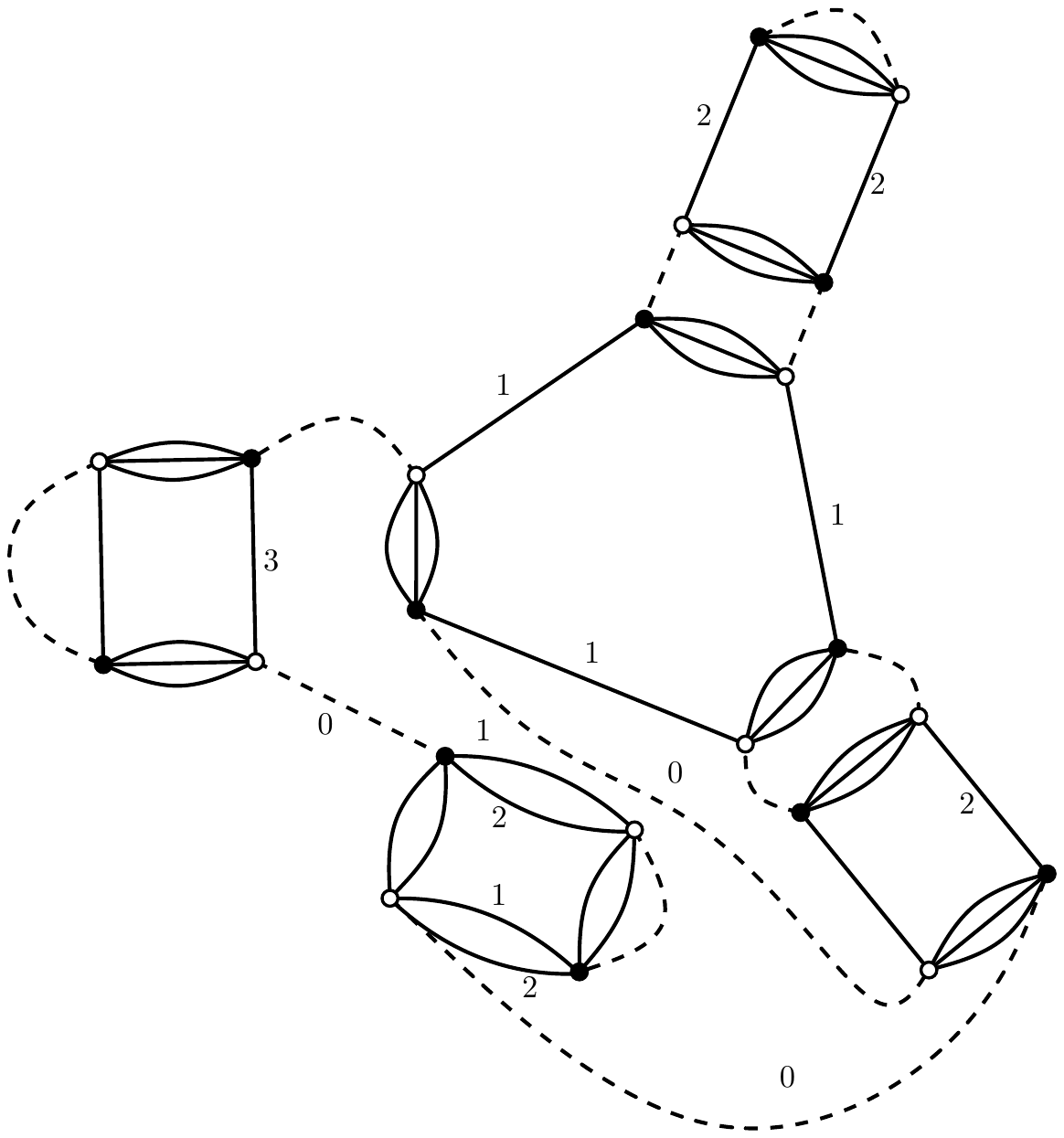}
\caption{\label{fig:VacFeynm} A closed 5-colored graph, where the dashed lines represent the edges of color 0.}
\end{figure}

If $\{\cB_\alpha\}_{\alpha\in A}$ is a finite set of bubbles, we denote $\bG_p(\{\cB_\alpha\}_{\alpha\in A})$ the set of $(D+1)$-colored graphs with colors in $\{0, \dotsc, D\}$ whose vertices have degree $D+1$ except for $p$ white vertices (and $p$ black vertices) which do not have an incident edge of color 0, and such that the connected components of the graphs obtained by removing all edges of color 0 are from the set $\{\cB_\alpha\}_{\alpha\in A}$ only. An example is given in the Figure \ref{fig:VacFeynm}.

A \emph{face of color $i\in [D]$} in a graph $\cG\in \bG_p(\{\cB_\alpha\}_{\alpha\in A})$ is a closed chain alternating edges of colors 0 and $i$. Vertices which do not have an incident edge of color 0 are called \emph{free vertices}. When $\cG$ has some free vertices, there are \emph{broken faces of color $i$} which are open chains of alternating colors 0 and $i$ between a white and a black free vertex. Broken faces can be naturally oriented from the white free vertex to the black free vertex.

\begin{figure}
\includegraphics[scale=.6]{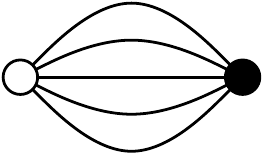} \hspace{2cm} \includegraphics[scale=.5]{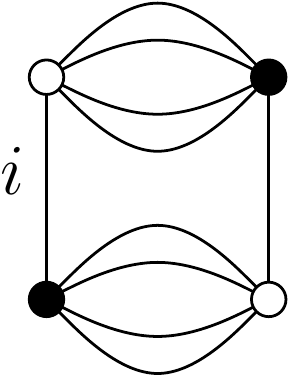}
\caption{\label{fig:BubbleExamples} There is a unique bubble with two vertices (on the left). On the right is the quartic melonic bubble $\cB_i$.}
\end{figure}

The quartic melonic bubble $\cB_i$ of color $i\in[D]$ has two white vertices $v_1, v_2$ and two black vertices $\bar{v}_1, \bar{v}_2$ such that $v_1$ is connected to $\bar{v}_1$ by $D-1$ edges (and similarly for $v_2$ and $\bar{v}_2$) of colors different from $i$ and $v_1$ is connected to $\bar{v}_2$ by an edge of color $i$, as well as $v_2$ and $\bar{v}_1$ (see Figure \ref{fig:BubbleExamples}).

\emph{Edge-colored maps} are combinatorial maps in which each edge carries a color in $[D]$. Note that for maps, we use edge-colorings which are typically non-proper. We denote $\bM_p$ the set of edge-colored maps with $p$ \emph{cilia}, i.e. marked corners, such that a vertex can have at most one cilium. A cilium divides a corner into two corners.

\begin{theorem} \label{thm:QuarticMelonicBijection}
There is a bijection between the sets $\bG_p(\{\cB_i\}_{i\in [D]})$ of $(D+1)$-colored graphs built from the quartic melonic bubbles $\{\cB_1, \dotsc, \cB_D\}$ and $\bM_p$.
\end{theorem}

We describe how this bijection works. First we recall how to describe a map $\cM = (\alpha, \sigma) \in\bM_p$ using two permutations $\alpha$ and $\sigma$. Label the half-edges of the map, including its cilia, of the map $1_{c_1}, 2_{c_2}, \dotsc, (2E)_{c_{2E}}$ where $c_k$ records the color of the half-edge $k$, and use primed labels to distinguish the cilia (which should be considered as having all the colors). As for ordinary combinatorial maps, the cycles of $\sigma$ are the cyclically ordered labels of the half-edges around each vertex. $\alpha$ acts as a fixed point free involution on the non-primed labels, which transposes labels of the same color only and thus describes the edges of the map, and has a fixed point at each cilium.

Let $\cG\in \bG_0(\{\cB_i\}_{i\in [D]})$. Each white vertex has a natural black partner: the one it is connected to via at least $D-1$ edges (it may be $D$ edges due to the color 0). Label each white vertex $1_{c_1}, \dotsc, (2E)_{c_{2E}}$ (and their black partner $\bar{1}_{c_1}, \dotsc, (\overline{2E})_{c_{2E}}$) where $c_k\in [D]$ is the color which does not connect the white vertex $k$ to its black partner $\bar{k}$. The permutation $\sigma$ encodes the edges of color 0: it is defined by $\sigma(a) = b$ if $a$ is a white vertex connected to $\bar{b}$ via an edge of color 0. One can thus think of $\sigma$ as mapping a white vertex $a$ to $\bar{b}$ following an edge of color 0, then to its white partner $b$. The cycles of $\sigma$ therefore represent cycles in $\cG$ made of edges of color 0 and parallel edges which connect partner vertices.

Furthermore, we get the involution $\alpha$ which respects the coloring by defining $\alpha(a)=b$ if $a$ labels a white vertex and $\bar{b}$ a black vertex which are not partners and if there is an edge of color $c\in [D]$ between them. In other words, $\alpha$ encodes the fact that the pairs of partners $(a, \bar{a})$ and $(\alpha(a), \overline{\alpha(a)})$ are connected by the color $c$ in each bubble $\cB_c$. 

\begin{figure}
\includegraphics[scale=.5]{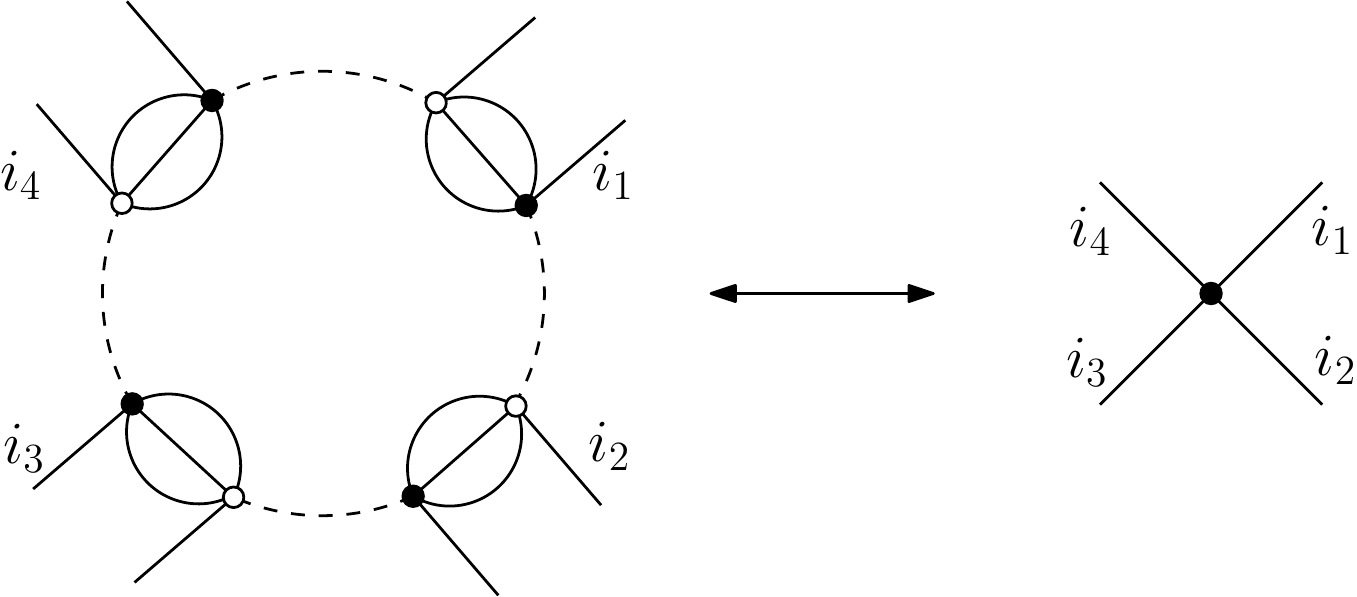}\hspace{2cm}
\includegraphics[scale=.5]{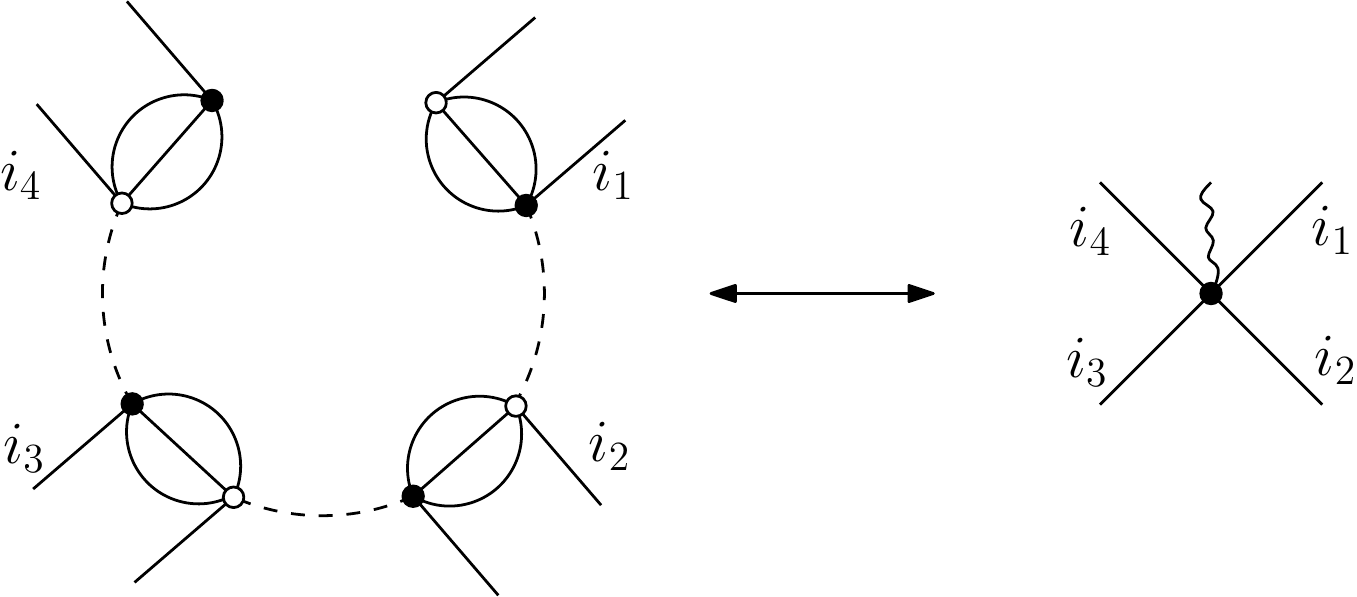}
\caption{\label{fig:IFVertex} Cycles which alternate edges of color 0 and parallel edges between partner vertices are represented in $\cM$ as vertices. Open chains become ciliated vertices where the position of the cilium indicates the missing edge of color 0.}
\end{figure}

A map $\cM\in\bM_0$ is then obtained from those two permutations $\alpha$ and $\sigma$. In other words, as illustrated in Figure \ref{fig:IFVertex},
\begin{itemize}
\item the melonic bubbles $\cB_i$ are mapped to edges of color $i$ of the map,
\item the cycles alternating edges of color 0 and $D-1$ parallel edges are mapped to vertices of the map,
\item which means that their corners represent the edges of color 0.
\end{itemize}

In the case where $\cG\in\bG_p(\{\cB_i\}_{i\in [D]})$ for $p\neq0$, the definition of $\sigma$ has to be amended. Start from a black vertex $\bar{v}$ with no incident edge of color 0. If its white partner $v$ also has no incident edge of color 0, we add an edge of color 0 between them, mark it and give it a primed label. This gives a cycle of $\sigma$ represented in $\cM$ as a vertex and the primed label is represented with a cilium. Else, $v$ has an incident edge of color 0 and we iterate $\sigma$ until we arrive at a white vertex with no incident edge of color 0. Then a marked edge of color 0 is added between this vertex and $\bar{v}$. This defines a cycle of $\sigma$ where the marked edge is promoted to a half-edge with a primed label. This cycle of $\sigma$ is a vertex of the map $\cM$ and the primed label becomes a cilium (see Figure \ref{fig:IFVertex}).

\begin{figure}
\includegraphics[scale=0.5]{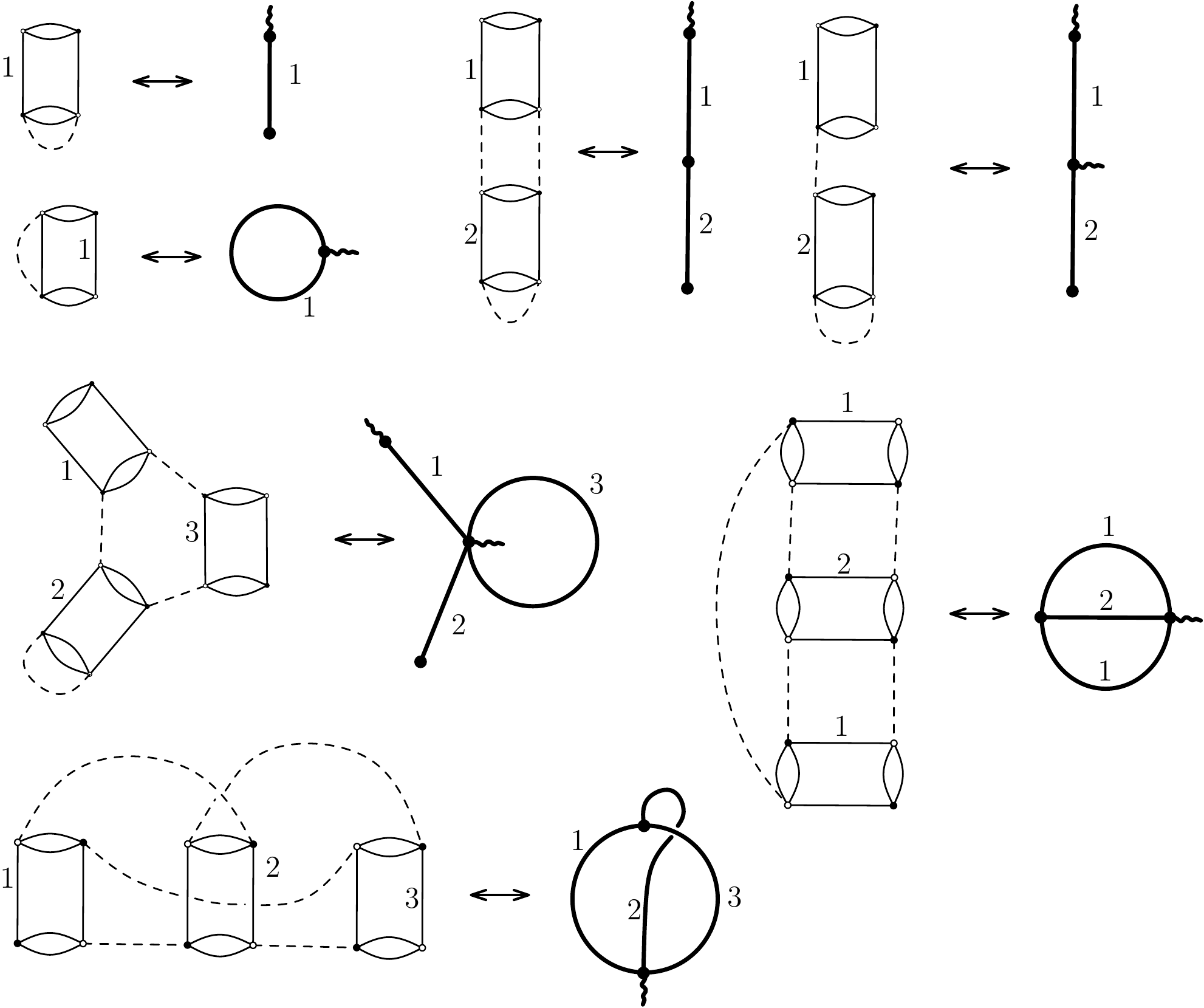}
\caption{\label{fig:EquMelEx} Here are some examples of the bijection between graphs with quartic melonic bubbles and edge-colored maps with cilia.}
\end{figure}

Examples of the bijection are given in Figure \ref{fig:EquMelEx}.

The \emph{monochromatic submap $\cM^{(i)}$} of color $i\in [D]$ is defined as the map obtained by removing all edges except those of color $i$. It is clear that the faces of color $i$ in $\cG$ become through the bijection the faces of $\cM^{(i)}$ (including its isolated vertices). If $\cM = (\alpha, \sigma)$, one defines $\alpha^{(i)}, \sigma^{(i)}$ by removing from the cycles of $\alpha$ and $\sigma$ all the labels of color different from $i$ (and keeping the primed labels, which correspond to cilia). Then the faces correspond to the cycles of $\alpha^{(i)} \sigma^{(i)}$. In the presence of cilia, a cycle of $\alpha^{(i)} \sigma^{(i)}$ may contain one or more primed labels. A broken face of color $i$ is then a part of a cycle between two cilia. It now receives its orientation from the cycle.

A melonic bubble is thus mapped to an edge of a map because it has two ``ends'', i.e. two pairs of partner vertices. Thus, for more generic bubbles, we expect ``edges'' with more than two ends, i.e. \emph{hyperedges}. For instance, it is clear that the bijection works similarly if instead of quartic melonic bubbles, we use melonic cycles with $n$ partner vertices, which are mapped to hyperedges with $n$ ends. Following Walsh \cite{Walsh}, hypermaps are in bijection with bipartite maps, and the melonic cycle of length $n$ can be represented as a white vertex of degree $n$ as shown in Figure \ref{fig:MonoColMel}. In this context, the maps which are obtained are called Walsh maps, and here they are edge-colored Walsh maps. However, for more generic bubbles, there is no preferred pairing of the vertices (they do not have natural partners). Moreover, the colored structure of a bubble cannot be represented as a simple white vertex {\it \`a la Walsh}, but instead we will show that a colored map can be used. Trading Walsh's vertices for maps leads to \emph{stuffed Walsh maps}.

\section{Bijection between edge-colored graphs and stuffed Walsh maps} \label{sec:Bijection}

\begin{definition}[Stuffed Walsh maps] \label{def:StuffedWalshMaps}
Let $\cM\in\bM_q$ be an edge-colored map with $q$ cilia. We denote $\bW_p(\cM)$ the set of stuffed Walsh maps of submaps $\cM$ and with $p$ cilia. A stuffed Walsh map $\cW\in \bW_p(\cM)$ is defined by
\begin{itemize}
\item black vertices which can only be connected to blue vertices,
\item removing all the black vertices and replacing their incident edges with cilia incident to blue vertices produces a disjoint union of edge-colored maps with cilia all isomorphic to $\cM$,
\item an edge which connects a black to a blue vertex carries a color subset of $[D]$ which is the set of colors incident to the blue vertex in $\cM$,
\item cilia are attached to black vertices, one cilium per vertex at most.
\end{itemize}
\end{definition}
Thus, while a black vertex is always connected to a blue vertex, blue vertices can be connected together. In fact, they form the submaps isomorphic to $\cM$. In the following, we will most of the time identify $\cM$ with the map where the cilia are replaced with half-edges, where a half-edge incident to a vertex of $\cM$ is labeled by a color set containing the colors incident to the vertex.

We also define the \emph{projected Walsh maps}, by contracting all the edges between blue vertices. This shrinks each submap $\cM$ to a point that can be represented as a white vertex and one obtains a bipartite map with white vertices of fixed degree and fixed incident color sets. If $\cW\in\bW_0(\cM)$, we denote $P(\cW)$ its unique projected map.

Consider a bubble $\cB$ with $\cV$ black vertices, and $\cG\in\bG_p(\cB)$ a $(D+1)$-colored graph whose maximal subgraph of colors $1, \dotsc, D$ is a disjoint union of bubbles $\cB$. We are going to build a bijection which maps $\cG$ to a stuffed Walsh map. This is done in three steps,
\begin{enumerate}
\item a universal part, independent of the details of $\cB$, which gives rise to the black vertices of the stuffed Walsh map,
\item a mapping of the bubble $\cB$ to an edge-colored map $\cM$ with $\cV$ cilia,
\item finally the gluing of those edge-colored maps to the universal black vertices.
\end{enumerate}

The bijection relies on a choice of pairing of $\cB$.
\begin{definition}[Pairings] \label{def:Pairing}
A \emph{pairing} of the vertices of a bipartite, edge-colored graph is a partition $\Omega$ of its vertices in pairs of black and white vertices.
\end{definition}

\subsection{The universal part of the bijection} \label{sec:BlackVertex}

Let $\Omega$ be a pairing of $\cB$. By following alternatively the edges between the paired vertices of each bubble in $\cG$ and the edges of color 0, we observe cycles which we represent as black vertices,
\begin{equation}
\begin{array}{c}\includegraphics[scale=.5]{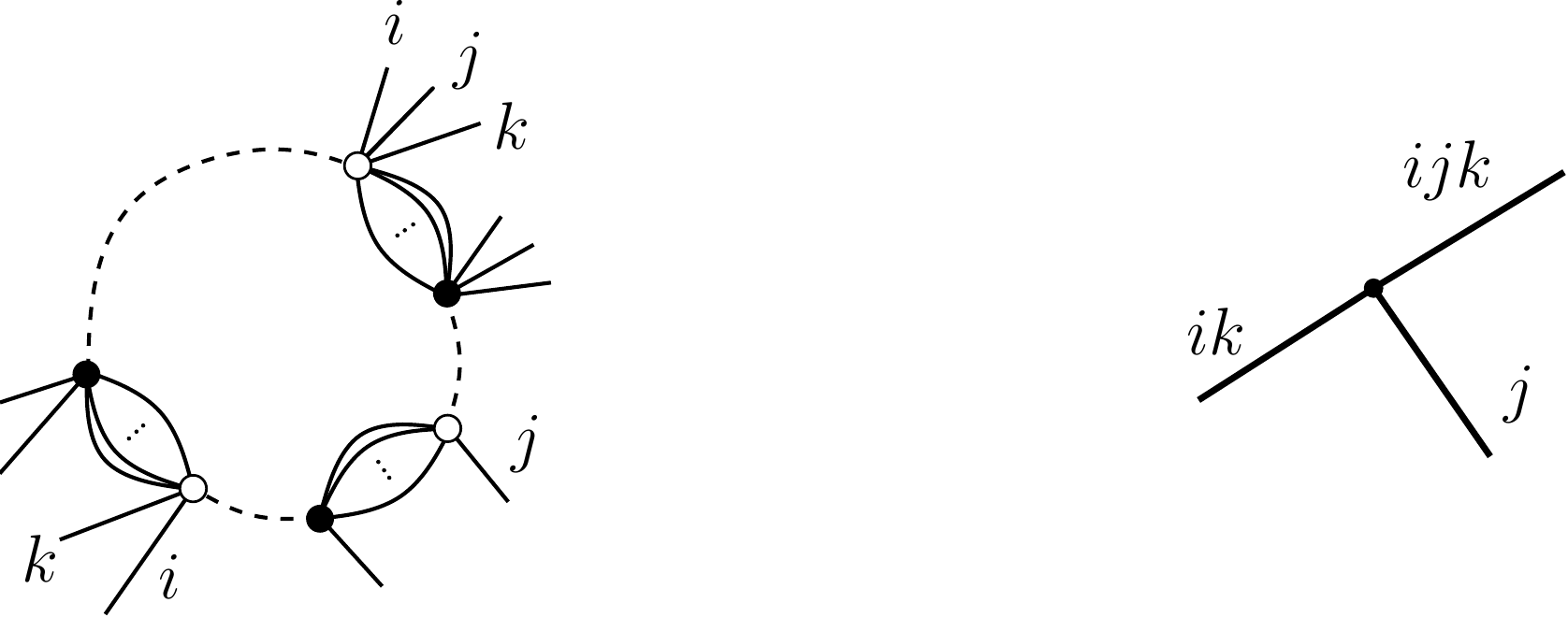}\end{array}
\end{equation}

More rigorously, those cycles can be defined as the cycles of a permutation. Denote $b$ the number of bubbles $\cB$ in $\cG$, and $\cV$ the number of white vertices of $\cB$.
\begin{itemize}
\item $\Omega$ induces a pairing $\Omega_\cG$ on the vertices of $\cG$: a white vertex $v$ of $\cG$ belongs to a single bubble $\cB$ in which it is paired to a black vertex $\bar{v}$ through $\Omega$, and they form the pair $(v, \bar{v})\in \Omega_\cG$.
\item An arbitrary labeling $1, \dotsc, b\cV$ of the white vertices of $\cG$ induces a labeling of the pairs of vertices of $\Omega_\cG$: the label $i$ goes to the pair whose white vertex has the label\footnote{One could use an arbitrary labeling of the black vertices and define a permutation $\tau_\Omega$ which maps the label of a white vertex to that of the black vertex it is paired with. Here we choose to set $\tau_\Omega = \operatorname{id}$.} $i$.
\item If $\cG$ is a closed graph, each vertex of $\cG$ is incident to an edge of color 0. Edges of color $0$ hence also induce a permutation $\mu$ on the labels $\{1,\dotsc,b\cV\}$: $x = \mu(y)$ if there is an edge of color 0 between the black vertex of the pair $y$ and the white vertex of the pair $x$.
\end{itemize}
This way, the cycles of $\cG$ which alternate edges of color 0 and parallel edges between paired vertices are precisely the cycles of the permutation $\mu$. We represent them as black vertices $\{v_\bullet\}$. Notice that this construction does not depend on the details of the bubble $\cB$. This is the reason why this part of the bijection is \emph{universal}.

Consider such a cycle represented as $v_\bullet$. At each iteration of $\mu$, arriving on a white vertex $i$ (or its black partner), we look at the colors of the edges incident to $i$ and which do \emph{not} connect it to its partner. We record those colors as the label $\cI\subset [D]$ of a half-edge incident to $v_\bullet$.

There is a cyclic order, inherited from the order of the pairs within the cycle of $\mu$ corresponding to $v_\bullet$, which we represent clockwise. This way, the corners around $v_\bullet$ correspond to the edges of color 0 along the cycle and they go clockwise from black to white vertices. Equivalently, a half-edge has two sides, and after following a corner clockwise, the side of the half-edge which is first met corresponds to a white vertex, while the other side of the half-edge corresponds to the black vertex of the same pair.

If $\cG$ is not a closed graph, i.e. $\cG\in\bG_p(\cB)$ for $p\neq 0$, then we proceed as follows. Consider a white vertex $v\in\cG$ with no incident edge of color 0. Either its black partner also has no incident edge of color 0, or it has one which defines the image of the pair by $\mu$. This continues until one arrives at a black vertex $\bar{v}'$ with no incident edge of color 0. We add an edge of color 0 between $v$ and $\bar{v}'$ and mark it. We thus get a cycle of $\mu$ which can be represented as a black vertex $v_\bullet$. Moreover, the marked edge of color 0 now corresponds to a corner of $v_\bullet$, which we mark with a cilium,
\begin{equation}
\begin{array}{c}\includegraphics[scale=.5]{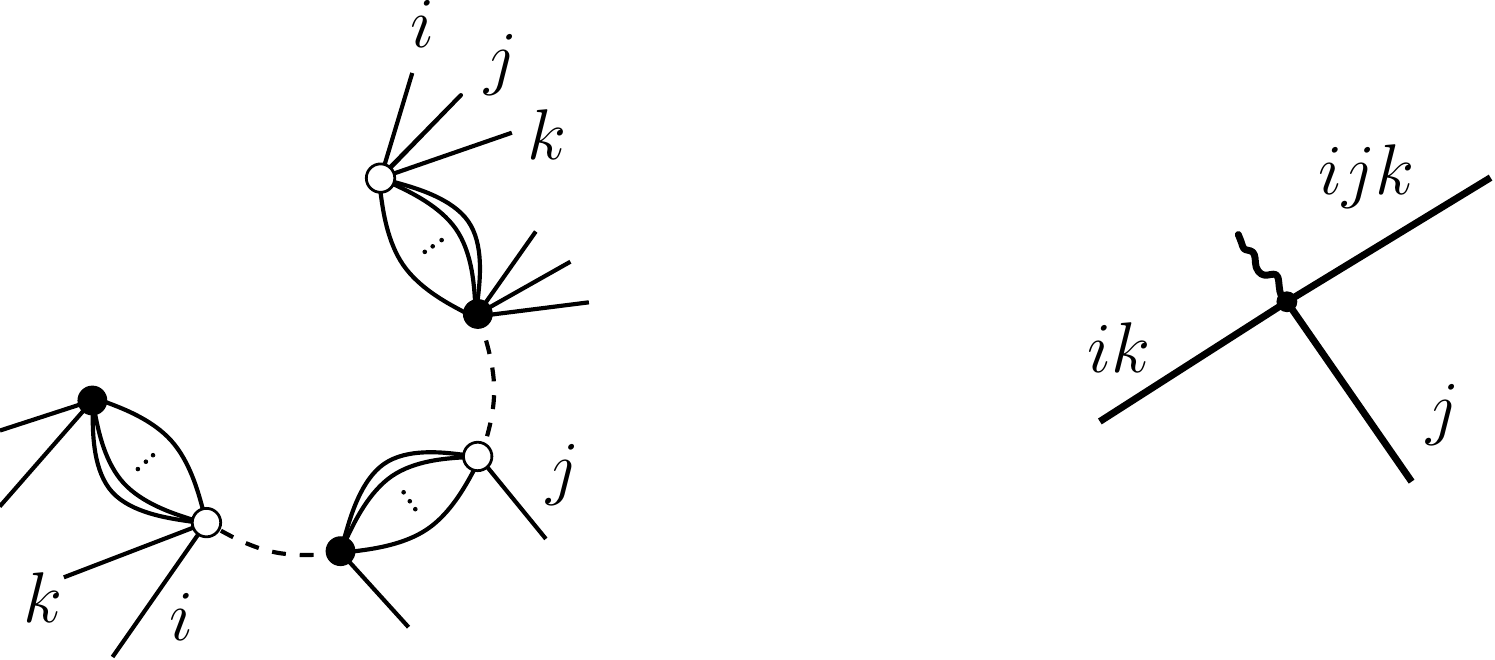}\end{array}
\end{equation}

\subsection{Bubbles as boundary graphs of edge-colored maps: the map $\cM(\cB, \Omega)$} \label{sec:BoundaryGraphs}

The black vertices $\{v_\bullet\}$ and their incident half-edges represent all pairs of vertices of $\cG$, all edges of color 0 and some of its cycles. There may be several cycles of $\mu$ going through the same bubble $\cB\subset \cG$. Those cycles are represented as black vertices $v_{\bullet 1}, v_{\bullet 2}, \dotsc$ which therefore have to be connected in some way.
\begin{itemize}
\item One can try and represent $\cB$ as a hyperedge connecting those vertices. Via Walsh's bijection \cite{Walsh} this hyperedge is just a white vertex whose incident half-edges correspond to the pairs of $\Omega$ and are labeled by the colors which do not connect the vertices in each pair. The order of the edges around the white vertex has to be fixed as they each represent a specific pair of vertices of $\cB$. By gluing the half-edges incident to black vertices to the half-edges \emph{of the same color labels} and incident to such white vertices we obtain \emph{projected Walsh maps}.

\item However, such a white vertex does not faithfully represent $\cB$. There can be several distinct bubbles for which the pairings are such that the incident edges have the same color labels. For instance, the two bubbles below, equipped with the pairing specified by the dashed ellipses,
\begin{equation}
\begin{array}{c} \includegraphics[scale=.5]{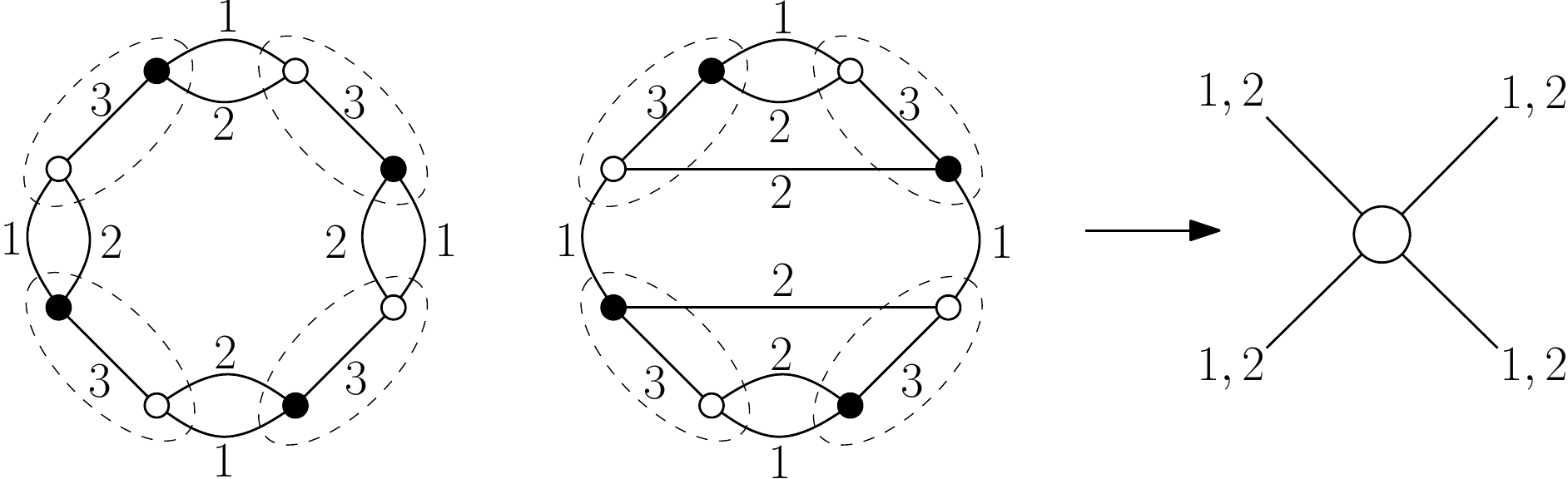}
\end{array}
\end{equation}
project onto the same white vertex. The projected Walsh maps coming from both $\bG_0(\cB), \bG_0(\cB')$ will thus be the same. This means that it is then impossible to reconstruct the graph $\cG$ without prior knowledge of the bubble $\cB$. For the same reason, it is impossible to track down the faces of color $i\in[D]$ through the bijection.

\item To remedy this issue, we will find a way to faithfully represent all faces of $\cG$. Indeed, $\cG$ is fully characterized by the knowledge of which colors connect which two vertices. Given the black vertices $\{v_\bullet\}$ that we have constructed and their incident half-edges, it therefore only remains to specify which colors (if any) connect two half-edges incident to black vertices. $\cB$ thus has to be represented as a \enquote{generalized, colored hyperedge} which dispatches the colors between its incident half-edges.
\end{itemize}

Since the black vertices $\{v_\bullet\}$ are locally embedded, we look for a representation of the bubble $\cB$ as a map. To do so, we will use the existing bijection for edge-colored graphs with quartic melonic bubbles.

\begin{definition}[Boundary graph]
The boundary graph $\partial \cG$ of an open Feynman graph $\cG\in\bG_p(\{\cB_\alpha\}_{\alpha\in A})$ has all the free vertices of $\cG$ as vertices, and an edge of color $i$ between two of them if there is a broken face of color $i$ (i.e. an open chain of colors $0$ and $i$) between them in $\cG$. 
\end{definition}

It is easy to see that the boundary graph $\partial \cG$ is a disjoint union of bubbles. Conversely, we show that all bubbles can be obtained as boundary graphs of edge-colored graphs with quartic melonic bubbles. Thanks to the bijection between edge-colored graphs with quartic melonic bubbles and edge-colored maps, it is then possible to represent any bubble $\cB$ as the boundary of a map. 

\begin{theorem} \label{thm:BoundaryGraphs}
The map $\partial: \bM_p \to \bB_p$ which takes the boundary graph of an edge-colored map with $p$ cilia is surjective on the set of bubbles with $p$ black vertices. Each cilium of an edge-colored map corresponds to a pair of vertices of the boundary bubble.
\end{theorem}

It is possible to prove the result for the boundary map $\partial: \bG_p(\{\cB_c\}_{c=1,\dotsc,D}) \to \bB_p$ and eventually use the bijection between $\bG_p(\{\cB_c\}_{c=1,\dotsc,D})$ and $\bM_p$. Instead we offer a proof where the boundary graphs are studied directly on $\bM_p$. In particular, for any bubble and pairing we construct in the proof a specific map $\cM(\cB, \Omega)$ such that
\begin{equation}
\partial \cM(\cB, \Omega) = \cB,
\end{equation}
and which will be used throughout the rest of the article.

{\bf Proof of Theorem \ref{thm:BoundaryGraphs}.} Consider a bubble $\cB$ and equip it with a pairing $\Omega$. We orient the edges of $\cB$ from white to black vertices. We now merge the vertices $(v, \bar{v})$ paired in $\Omega$ into new vertices $V=(v, \bar{v})$ and erase the edges which connect the two vertices of each pairs, to obtain the graph $\cB_{\circlearrowleft, \Omega}$. Clearly, $\cB_{\circlearrowleft, \Omega}$ is equivalent to $\cB$.

Both $\cB$ and $\cB_{\circlearrowleft, \Omega}$ are characterized by a set of permutations $\tau_1, \dotsc, \tau_D$. Label the pairs of vertices of $\Omega$ and define $\tau_i$ on $\cB$ as the map which sends the pair of vertices $x$ to the pair $y$ if there is an edge of color $i$ between the white vertex of the pair $x$ and the black vertex of the pair $y$. On $\cB_{\circlearrowleft, \Omega}$, $\tau_i$ maps $x$ to $y$ if there is an edge of color $i$ directed from $x$ to $y$.

The subgraph $\cB_{\circlearrowleft, \Omega}^{(i)}$ obtained by removing all edges but those of color $i$ has vertices of degree zero or two, and is thus a disjoint union of directed cycles $C_1^{(i)}, \dotsc, C_K^{(i)}$, each of them corresponding to a cycle of the permutation $\tau_i$.

To eventually represent $\cB$ as a map, we have to specify a way of embedding the graph $\cB_{\circlearrowleft, \Omega}$ and marking corners while erasing the orientations. Then we will see that the embedding determines the permutations $\alpha, \sigma$ which formally define the map. 

\begin{figure}
\includegraphics[scale=0.9]{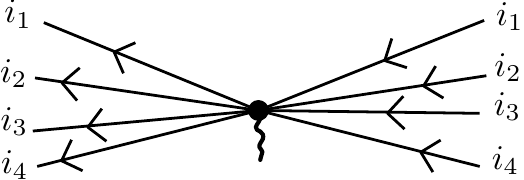}
\caption{\label{fig:FlowVert} Going clockwise around the vertex, one encounters the cilium, then all the outgoing edges and then all the ingoing edges.  }
\end{figure}

There has to be exactly one cilium on each vertex. The key idea is to embed each cycle $C_k^{(i)} \mapsto \tilde{C}_k^{(i)}$, for $k=1, \dotsc, K$, such that all cilia attached to $\tilde{C}_k^{(i)}$ lie on a single face. To do so, we split the local neighborhood around each vertex into two regions, and draw the ingoing edges in one region, and the outgoing edges in the other region. The relative order between edges of distinct colors can be chosen arbitrarily (we will see it does not affect the boundary graph). A cilium is finally added between the two regions, such that the outgoing edges are encountered first when one goes around the vertex clockwise starting from the cilium (see Figure \ref{fig:FlowVert}). This prescription ensures that the orientations can be safely erased. Indeed, the first time a color is encountered when one goes clockwise around the vertex starting from the cilium, we know it is outgoing. Moreover, this allows to identify the side of the cilium adjacent to the outgoing edges as corresponding to the white vertex of the pair, and the other side of the cilium to the black vertex.

\begin{figure}
\includegraphics[scale=0.55]{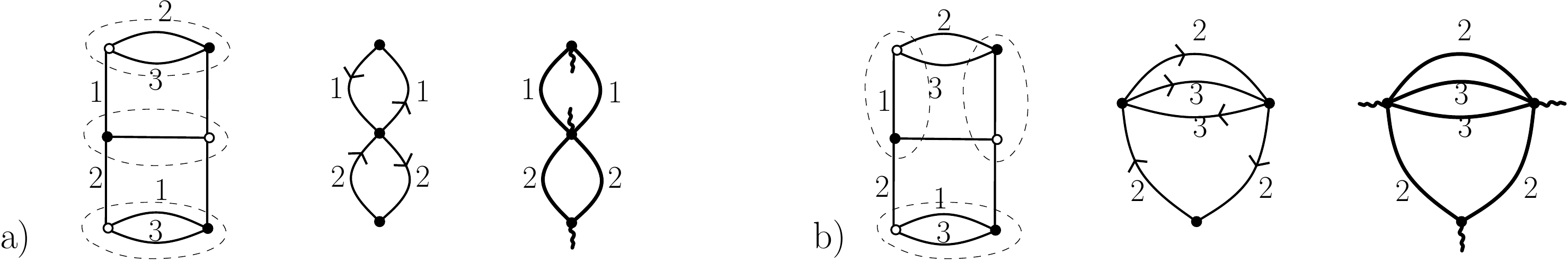}
\caption{\label{fig:FromBubbleToMap} In both cases $a)$ and $b)$, we have represented from left to right the bubble $\cB$ and its pairing, the graph $\cB_{\circlearrowleft, \Omega}$ obtained by contracting the edges between the paired vertices and orienting the edges from white to black, and then its embedding as a ciliated map, so that all cilia lie on the same face of each cycle for every monochromatic submap $\cM^{(i)}$.}
\end{figure}

We thus obtain an edge-colored map $\cM(\cB, \Omega)$. An example of the construction is given in Figure \ref{fig:FromBubbleToMap}, with two different pairings for the same bubble. As in Section \ref{sec:Review}, $\cM(\cB, \Omega)$ has to be characterized by a pair of permutations, $\cM(\cB, \Omega) = (\alpha, \sigma)$. Since we have ignored the relative order between edges of different colors around each vertex, it is sufficient to identify the permutations $\alpha^{(i)}, \sigma^{(i)}$ for $i=1, \dotsc, D$, which are the restrictions of $\alpha$ and $\sigma$ to the half-edges of color $i$ (including the cilia).

For a labeling $\{x\}$ of the pairs of vertices of $\cB$, on which the permutations $\tau_i$ are defined, we define a labeling of the half-edges of $\cB_{\circlearrowleft, \Omega}$: $h^{(i)}_{a,x}$ is the half-edge of color $i$ incident to $x$, with $a=1$ if it is outgoing at $x$ and $a=2$ if it is incoming.
\begin{itemize}
\item The permutation $\alpha^{(i)}$ describes the gluing of two half-edges into an edge. Since the edges of $\cM(\cB, \Omega)$ are the same as those of $\cB_{\circlearrowleft, \Omega}$, we define 
\begin{equation}
\alpha^{(i)}(h^{(i)}_{1,x}) = h^{(i)}_{2,y} \qquad \text{and} \qquad \alpha^{(i)}(h^{(i)}_{2,y}) = h^{(i)}_{1,x} \qquad \text{if $\tau_i(x) = y$}.
\end{equation}
We supplement it with a fixed point at each cilium.
\item As for $\sigma^{(i)}$, it encodes the order of the edges of color $i$ and the cilium at each vertex. Since at each vertex, $\cB_{\circlearrowleft, \Omega}$ is either of degree zero or two on the color $i$, so is the monochromatic submap $\cM^{(i)}$. Thus, either it gives a trivial cycle (corresponding to the cilium) or a cycle of length three between the cilium and the two incident half-edges. The cyclic order between the three of them has been specified by the embedding of the cycles $C^{(i)}_k$. It gives
\begin{equation}
\sigma^{(i)} = \prod_x (\text{cilium}_x, h^{(i)}_{1,x}, h^{(i)}_{2,x}).
\end{equation}
\item Furthermore, the faces of color $i$ are given by the cycles of $\alpha^{(i)}\sigma^{(i)}$. For a vertex with just a cilium in $\cM^{(i)}$, one gets a cycle of length one. For vertices with two incident edges of color $i$, the embedding of the incident cycle $\tilde{C}^{(i)}_k$ gives rise to two faces: one has no cilia while the other has a cilium at each vertex. Indeed, if one starts with a cilium at $x_0$, one gets a cycle of the type
\begin{equation} \label{CycleBrokenFaces}
(\text{cilium}_{x_0}, h^{(i)}_{1,x_0}, h^{(i)}_{2,x_1}, \text{cilium}_{x_1}, h^{(i)}_{1,x_1}, h^{(i)}_{2, x_2}, \text{cilium}_{x_2}, \dotsc, \text{cilium}_{x_q}, h^{(i)}_{1,x_q}, h^{(i)}_{2,x_{q+1}}, \text{cilium}_{x_{q+1}}, \dotsc ),
\end{equation}
where the vertex $x_{q+1}$ is the vertex which comes after $x_q$ along the cycle, i.e. $\tau_i(x_q) = x_{q+1}$. If one starts a cycle from the half-edge $h^{(i)}_{1,x_0}$, no cilium is encountered: $(h^{(i)}_{1,x_0}, h^{(i)}_{2,x_0}, h^{(i)}_{1,x_1}, h^{(i)}_{2,x_1}, \dotsc, h^{(i)}_{1,x_q}, h^{(i)}_{2,x_q}, h^{(i)}_{1,x_{q+1}}, \dotsc)$.
\end{itemize}

We can now verify that $\partial \cM(\cB, \Omega) = \cB$. The boundary graph is obtained by drawing a pair of white and black vertices for each vertex of the map and drawing an edge of color $i$ from a white vertex $v$ to a black vertex $\bar{v}$ if there is a positively oriented broken face of color $i$ between the vertices of $\cM(\cB, \Omega)$ from which $v$ and $\bar{v}$ come.

As seen in Section \ref{sec:Review}, a broken face of color $i$ in an edge-colored map $\cM$ is a part of a cycle of $\alpha^{(i)} \sigma^{(i)}$ between two cilia. Then, through each vertex, there is at most one non-trivial cycle of $\alpha^{(i)} \sigma^{(i)}$ with cilia, coming from a cycle $C^{(i)}_k$ of $\cB_{\circlearrowleft, \Omega}$ and given in \eqref{CycleBrokenFaces}. Keeping only the ordered vertices, it reads $(x_1, x_2, \dotsc, x_{|C^{(i)}_k|})$ where $\{x_q\}$ labels the vertices of the cycle. That information is in fact sufficient to reconstruct $\cB$ since from $x_q$ to $x_{q+1}$ one goes from a white to a black vertex (this is equivalent to the labeling by $a=1, 2$ of the half-edges). Then we define the permutation $\tilde{\tau}_i$ by $\tilde{\tau}_i(x_q) = x_{q+1}$ (cyclically with $q$). It is well-defined since there is at most one non-trivial cycle of color $i$ through each vertex. Finally, we observe that
\begin{equation} \label{BrokenFacesVsEdges}
\tilde{\tau}_i = \tau_i,
\end{equation}
where $\tau_i$ is the permutation of color $i$ which characterizes the edges of color $i$ of $\cB$. This shows that $\partial \cM(\cB, \Omega) = \cB$.
\qed

In the above proof we have constructed a map $\cM(\cB, \Omega)$. There is however an infinite number of maps $\cM$ such that $\partial \cM = \cB$ for a fixed $\cB$ (add vertices of degree two to the edges of $\cM(\cB, \Omega)$ for instance). We provide below a graphical procedure which extracts the boundary graph of any map $\cM\in\bM_p$. The procedure obviously relies on broken faces, introduced in Section \ref{sec:Review}. Recall that a broken face is a walk of fixed color from one cilium to another (possibly the same). A broken face can be oriented by following the walk clockwise around each vertex of the map. Notice that there exist \emph{trivial} broken faces, i.e. which go around a single cilium without following an edge. The others are \emph{non-trivial} broken faces.

\begin{figure}
\includegraphics[scale=0.6]{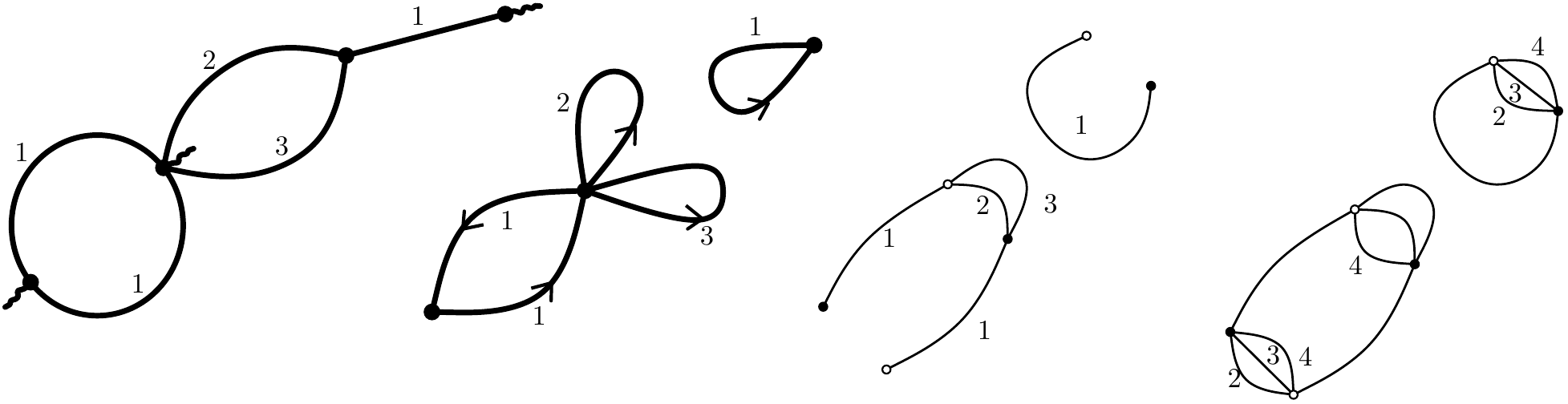}
\caption{\label{fig:boundM} The steps of the procedure giving the boundary graph of a map in $\mathbb{M}_p$. On the left is the map $\cM$ and to its right the graph $\cG_{\circlearrowleft, \cM}$. The third graph is obtained from $\cG_{\circlearrowleft, \cM}$ by replacing each vertex with a pair of black and white vertices. Adding the missing colors between the vertices of each pair leads to $\cB = \partial \cM$.}
\end{figure}

The following procedure extracts the boundary graph of an edge-colored map.
\begin{enumerate}
\item Keep only the ciliated vertices in $\cM$. For each non-trivial broken face of color $i$ between two ciliated vertices, draw an edge of color $i$ between them. Each of them is oriented as its corresponding broken face. One thus gets a directed graph $\cG_{\circlearrowleft, \cM}$.

Recall from Section \ref{sec:Review} that if $\cM = (\alpha, \sigma)$, the faces of color $i$ are the cycles of $\alpha^{(i)}\sigma^{(i)}$ and the broken faces are the parts of those cycles between two cilia. This shows that the subgraph of $\cG_{\circlearrowleft, \cM}$ containing only edges of color $i$ consists of disjoint oriented cycles.
\item The next step is to replace each vertex $V$ by a pair of white and black vertex $(v, \bar{v})$, so that all the incoming (respectively outgoing) edges attached to $V$ now reach $\bar{v}$ (resp. $v$). At this stage, edge orientations may be erased. 
\item For each pair $(v, \bar{v})$ coming from a vertex $V$ of $\cM$, edges are added if necessary between $\bar{v}$ and $v$ so that each vertex has all $D$ colors incident. The corresponding non-directed graph is the boundary graph of $\cM$. 
\end{enumerate}
The steps of this procedure are described in Figure \ref{fig:boundM}. Note that the two last steps induce a pairing $\Omega_\cM$ of $\partial \cM$: the paired vertices are those which come from the same vertex of $\cM$ (or $\cG_{\circlearrowleft, \cM}$). We call it the \emph{pairing induced by $\cM$}.

The above constructions have introduced two types of directed graphs: $\cB_{\circlearrowleft, \Omega}$ obtained from a bubble $\cB$ equipped with a pairing $\Omega$, and $\cG_{\circlearrowleft, \cM}$ obtained from a map $\cM$. In the proof of Theorem \ref{thm:BoundaryGraphs}, we have built a map $\cM(\cB, \Omega)$ of boundary graph $\cB$. In fact, the equality \eqref{BrokenFacesVsEdges} $\tilde{\tau}_i = \tau_i$ for each color $i\in [D]$ is equivalent to 
\begin{equation}
\cG_{\circlearrowleft, \cM(\cB, \Omega)} = \cB_{\circlearrowleft, \Omega}.
\end{equation}
More generally, if one only assumes $\cB = \partial \cM$, then $\cB_{\circlearrowleft, \Omega}$ and $\cG_{\circlearrowleft, \cM}$ have the same number of vertices but may be different. Yet, if the pairing $\Omega_\cM$ induced by $\cM$ is the same as the one equipping $\cB$, $\Omega_\cM = \Omega$, then the two graphs can only differ by the presence of loop-edges in $\cG_{\circlearrowleft, \cM}$ (which by construction $\cB_{\circlearrowleft, \Omega}$ cannot have).

The set of graphs $\{\cB_{\circlearrowleft, \Omega}\}$ obtained from bubbles and pairings consists of non-properly-edge-colored, directed, connected (since bubbles are connected) Eulerian graphs. Furthermore, they have vertices of even degrees less than or equal to $2D$ such that for each edge of color $i$ going in a vertex, an edge of the same color comes out of the same vertex, and each color cannot be incident more than twice on a vertex. Moreover, $\cB_{\circlearrowleft, \Omega}$ cannot have loop-edges, since the edges between paired vertices have been contracted down to a vertex. A graph from this set defines a unique bubble and pairing. 

\subsection{Stuffed Walsh maps}

\begin{theorem}
For any bubble $\cB$ equipped with a pairing $\Omega$, there is a bijection between $\bG_p(\cB)$ and $\bW_p(\cM(\cB, \Omega))$. The bubbles isomorphic to $\cB$ are mapped to the submaps isomorphic to $\cM(\cB, \Omega)$, the edges of color 0 to the corners around black vertices, and the faces of each color to the faces of the same color which go along at least one black vertex.
\end{theorem}

In fact, any map $\cW\in \bW_p(\cM)$ for any edge-colored map $\cM$ such that $\partial \cM = \cB$ represents a unique graph $\cG\in \bG_p(\cB)$. Also any graph $\cG\in \bG_p(\cB)$ can be represented as a stuffed Walsh map for any given choice of $\cM$ such that $\cB = \partial \cM$. However, the stuffed Walsh map may not be unique if $\cM$ does not have the symmetries of $\cB$ under the exchange of pairs of vertices. This is the reason why we use the specific map $\cM(\cB, \Omega)$.

{\bf Proof.} We prove the bijection explicitly for $p=0$. Let $\cG\in\bG_0(\cB)$ with $b$ bubbles and let $\cV$ be the number of white vertices of $\cB$. We label the bubbles in $\cG$ as $\kappa = 1, \dotsc, b$ and the pairs of vertices of $\cB$ in $\Omega$ as $x = 1, \dotsc, \cV$. Since the $b$ bubbles of $\cG$ are isomorphic to $\cB$, we can find an isomorphism as graphs with labeled pairs, thus inducing a labeling $1_\kappa, \dotsc, \cV_\kappa$ of the pairs of the bubble $\kappa$. 

We have seen in Section \ref{sec:BlackVertex} that each pair of vertices appears in a cycle made of edges of color 0 and parallel edges. Each such cycle is represented as a black vertex whose incident half-edges represent the pairs of vertices and whose corners represent the edges of color 0. This way, we represent all pairs as labeled half-edges incident on black vertices.

Each bubble $\kappa$ with labeled pairs of vertices is represented as a map $\cM(\cB, \Omega)_\kappa$. We color its vertices in blue to distinguish them from the black vertices, as in Definition \ref{def:StuffedWalshMaps}. Its cilia receive the labels of the pairs they correspond to. Each cilium can be replaced with a half-edge carrying all the colors incident to the vertex it is attached to.

Therefore each pair of vertices $x_\kappa$, for $x=1, \dotsc, \cV$ and $\kappa = 1, \dotsc, b$ is represented both as a half-edge incident to a black vertex and as a half-edge incident to a blue vertex of $\cM(\cB, \Omega)_\kappa$. We glue the half-edges which correspond to the same pairs (without twist), and this way obtain a stuffed Walsh map $\cW(\cG)$, where the relabeling of the pairs in $\cG$ correspond to relabelings of the edges between black and blue vertices of $\cW(\cG)$.

This construction may however seem ambiguous, which could result in the fact that a graph $\cG$ could be represented by several distinct stuffed Walsh maps. The only apparent ambiguity lies in the isomorphism between the bubbles $\kappa$ and $\cB$ as graphs with labeled pairs. Given one such isomorphism, they may be others which differ by some permutations $S$ on the labels $1, \dotsc, \cV$. Such symmetries of $\cB$ still remain symmetries of $\cM(\cB, \Omega)$ by construction, i.e. $\cM(\cB, \Omega)$ is invariant under the action of $S$ on the labeled cilia (or half-edges). It means that there is a unique way to glue each submap $\cM(\cB, \Omega)_\kappa$ to the black vertices. Thus all the possible isomorphisms lead to the same stuffed Walsh map and the construction is uniquely defined.

Conversely, each stuffed Walsh map $\cW\in\bW_0(\cM)$ where $\partial \cM = \cB$ gives rise to a colored graph $\cG\in \bG_0(\cB)$. It is obtained by taking the boundary of $\cM$ to extract $\cB$, and use the Section \ref{sec:BlackVertex} to reconstruct the edges of color 0 of $\cG$ out of the black vertices of $\cW$.

For $p\neq 0$, the bijection works similarly, using the part of Section \ref{sec:BlackVertex} on open cycles which map to black vertices with cilia.

By construction, the bijection maps the bubbles $\kappa = 1, \dotsc, b$ of $\cG\in \bG_p(\cB)$ to the submaps $\cM(\cB, \Omega)_\kappa$ in $\cW\in\bW_p(\cM(\cB, \Omega))$, and the edges of color 0 of $\cG$ to the corners around the black vertices of $\cW$. We will study the faces below in Section \ref{sec:Faces}.
\qed

Quite clearly, the bijection extends to colored graphs built from several distinct bubbles. Each of them can be represented as a map with cilia. They are then assembled via black vertices as constructed in Section \ref{sec:BlackVertex}, whose presence is universal.

\subsection{Face preserving and the reduced map $\tilde{\cM}(\cB, \Omega)$} \label{sec:Faces}

Using a slight modification to $\cM(\cB, \Omega)$, the faces of each color can also be preserved. First we have to show that the faces of $\cG$ become the faces of the corresponding stuffed Walsh map which go along at least one black vertex. To do so, we extend the definition of the monochromatic submap of an edge-colored map to stuffed Walsh maps: $\cW^{(i)}$, for $i\in[D]$, is the (typically disconnected) submap obtained from $\cW$ by removing all the edges whose color labels do not contain $i$. $\cW^{(i)}$ is an ordinary combinatorial map, and its faces are the faces of color $i$ of $\cW$.

Let $\cG\in\bG_0(\cB)$ and $\cW$ the corresponding stuffed Walsh map. The faces of color $i$ in $\cG$ become the faces of color $i$ in $\cW$ which go along at least one black vertex. Indeed, a face of $\cG$ necessarily goes along an edge of color 0, which becomes a corner of a black vertex in $\cW$. After such a corner, we meet a half-edge which corresponds to a pair of vertices of $\cB$. More precisely, by following the corner clockwise, the side of the half-edge which is met after the corner corresponds to the white vertex of the pair (see Section \ref{sec:BlackVertex}), while the other side of the half-edge corresponds to the black vertex. We have to check that the construction glues to that half-edge a half-edge attached to a blue vertex such that the black and white vertices are correctly identified. It can be checked from Section \ref{sec:BoundaryGraphs} that when following a corner of $\cM(\cB, \Omega)^{(i)}$ in a counter-clockwise fashion, the side of the cilium which is met after the corner corresponds to a white vertex of $\cB$. Therefore, by gluing the two half-edges without twist we correctly identify the two copies of the white and black vertices of each pair and thus reconstitute the faces of $\cG$ which go through this pair of vertices.

However, if $\cM(\cB, \Omega)$ contains (non-broken) faces of color $i$, then they give rise to faces of $\cW$ only bounded by blue vertices and which are not present in $\cG$. Therefore, we will modify $\cM(\cB, \Omega)$ so that the monochromatic submap $\cM^{(i)}(\cB, \Omega)$ becomes a forest. This can be done in one of two possible ways. Recall that the monochromatic submap $\cM^{(i)}(\cB, \Omega)$ consists of a disjoint union of cycles, where each cycle is incident to a face with no cilia. These faces are the supernumerary faces of $\cW$. To get rid of them, we can simply remove an edge of color $i$ from each cycle. That merges the face with no cilia with the other face incident to the cycle which carry all the cilia.

Another way to remove the supernumerary faces is to add a vertex for each cycle of $\cM^{(i)}(\cB, \Omega)$ and connect it to the vertices of the cycle while removing the edges of the cycle. The order of the edges around the new vertex is the order of the vertices along the cycle. This turns each cycle into a star and the monochromatic submap $\cM^{(i)}(\cB, \Omega)$ becomes a forest of stars.

Thus, both methods lead to a map $\tilde{\cM}(\cB, \Omega)$ which we call a \emph{reduced map} for $(\cB, \Omega)$, whose monochromatic submaps are forests. This means that the faces of $\cG$ can be mapped to the faces of $\cW\in\bW_p(\tilde{\cM}(\cB, \Omega))$.

\begin{figure}
\includegraphics[scale=.7]{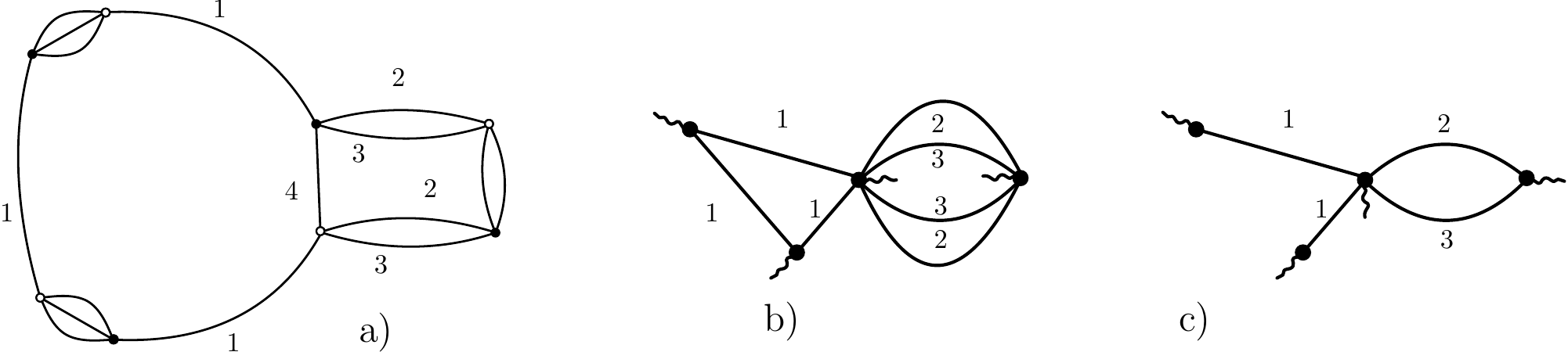}
\caption{\label{fig:FromBtoV} $b)$ is the map $\cM(\cB, \Omega)$ for the bubble in $a)$. In $c)$ we find a reduced map $\tilde{\cM}(\cB, \Omega)$ obtained by removing an edge of color 1, an edge of color 2 and an edge of color 3.}
\end{figure}

\begin{figure}
\includegraphics[scale=.6]{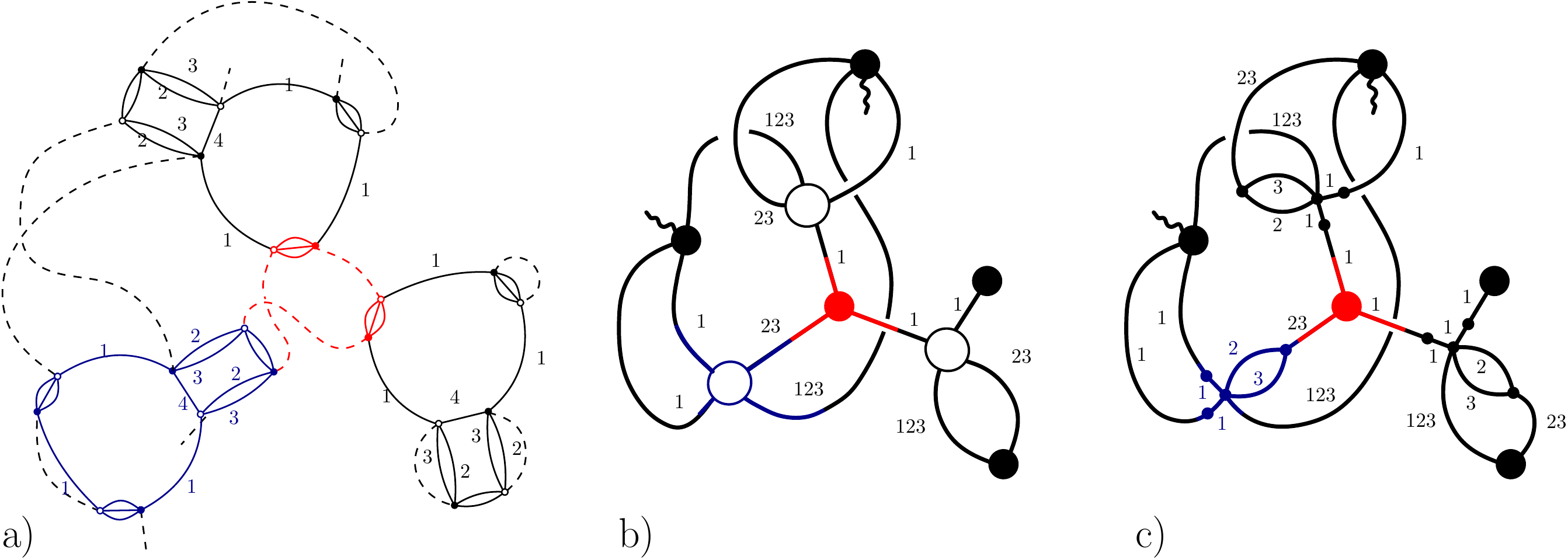}
\caption{\label{fig:Bijection1_2} $a)$ is an example of a colored graph with three bubbles isomorphic to the one of Figure \ref{fig:FromBtoV}. It gives rise to a map $\cW\in\bW_2(\tilde{\cM}(\cB, \Omega))$ in $c)$ where we have used the reduced map $\tilde{\cM}(\cB, \Omega)$ of Figure \ref{fig:FromBtoV}. In $b)$, this is the projected version of the map.}
\end{figure}

An example of a bubble and one of its reduced maps is given in Figure \ref{fig:FromBtoV}. A colored graph using the same bubble together with a stuffed Walsh map using the corresponding reduced map are shown in Figure \ref{fig:Bijection1_2}.

\subsection{Recovering the bijection between edge-colored graphs with quartic melonic bubbles and edge-colored maps}

We now compare the bijection we have constructed in the case of edge-colored graphs with quartic melonic bubbles with the known bijection presented in Section \ref{sec:Review}. The quartic melonic bubbles of color $i=1, \dotsc, D$ have two possible pairings: one denoted $\Omega$ which pairs a white vertex to the black vertex it is connected to by $D-1$ parallel edges of colors different from $i$, and one denoted $\Omega'$ which pairs it to the black vertex it is connected to by an edge of color $i$. They lead to two different bijections.

In the first case, the map $\cM(\cB_i, \Omega)$ is a cycle of length two with two cilia on the same face. Therefore $\tilde{\cM}(\cB_i, \Omega)$, obtained by removing one edge and replacing the cilia by half-edges of color $i$, consists of an edge with two half-edges attached to two bivalent blue vertices. It can obviously be simplified to a single edge of color $i$ and therefore stuffed Walsh maps reduce to edge-colored maps (there are no blue vertices and all the edges carry a single color). This is the bijection of Section \ref{sec:Review}.

The map $\tilde{\cM}(\cB_i, \Omega')$ also consists of an edge with two cilia, but the edge now gets the complementary label $\cI = [D]\setminus \{i\}$. This set of maps is obviously the same as with the pairing $\Omega$, except for the color labels.

The difference between the two pairings appears upon investigating the number of faces of those maps. For instance, the graph $\cG$ obtained by adding edges of color 0 between the vertices paired in $\Omega$ has $2(D-1) + 1$ faces. It is represented in $\bW_0(\tilde{\cM}(\cB_i, \Omega))$ as a map $m$ with an edge joining two black vertices. Since this edge has the color $i$, the monochromatic submap of color $i$ has a single face around the edge, while the $(D-1)$ other monochromatic submaps of color $j\neq i$ consists of two isolated vertices, hence reproducing the counting of faces of $\cG$. $\cG$ is also represented in $\bW_0(\tilde{\cM}(\cB_i, \Omega'))$ as a loop edge with color set $\cI$ and a single vertex. On the other hand, the map $m$, seen in $\bW_0(\tilde{\cM}(\cB_i, \Omega'))$, has a single non-loop edge with color set $\cI$, hence it has $(D-1)$ monochromatic submaps with a single face. Only the monochromatic submap of color $i$ has two isolated vertices and two faces, leading to an overall number of $2+(D-1)$ faces. Therefore the loop edge has more faces than the non-loop edge with $\Omega'$, while it works the other way around with $\Omega$.

\section{Dominant maps} \label{sec:DominantMaps}

Planar combinatorial maps are those which maximize the number of faces at fixed numbers of edges and vertices, or for regular maps, e.g. triangulations, quadrangulations, they maximize the number of faces at fixed number of vertices. This idea can be extended to colored graphs. Among all colored graphs, the maximal number of faces for a fixed number of white vertices is reached for melonic graphs, and only them \cite{Melons}. If $\cB$ is a melonic bubble with $p$ white vertices and $\cG\in\bG_0(\cB)$ a melonic graph with $b$ bubbles, then the number of faces of $\cG$ is
\begin{equation} \label{Faces}
\cF(\cG) \equiv \sum_{i=1}^D \cF^{(i)}(\cG) = (D-1)(p-1)b + D.
\end{equation}
The quantity $\omega(\cG) = D + (D-1)(p-1)b - \cF(\cG)$, called the reduced degree, is a positive integer which vanishes if and only if $\cG$ is melonic. However, if $\cB$ is not melonic, the set $\bG_0(\cB)$ does not contain melonic graphs.

We call the \emph{dominant graphs} in $\bG_0(\cB)$ those which maximize the number of faces at fixed number of bubbles (or vertices, since the number of white vertices is $pb$). Since the bijection preserves the faces, we can instead look for the dominant stuffed Walsh maps. In fact, projected Walsh maps will be the most useful representation in this section, so that we will call \emph{dominant maps} the dominant projected Walsh maps.

\subsection{Choice of pairing}

As seen in Section \ref{sec:Review} with quartic melonic graphs, the structure of dominant maps depends on the pairing chosen for the bubble $\cB$. In the quartic melonic case, the pairing $\Omega$ is such that the dominant maps are exactly the trees. However, if one uses $\Omega'$ instead, dominant maps will look different and have cycles.

Since our bijection extends the existing bijection for edge-colored graphs with quartic melonic bubbles, it seems reasonable to look for a pairing which makes trees dominant maps (notice that in general the stuffed Walsh maps cannot be trees since the submap $\tilde{\cM}(\cB, \Omega)$ itself is generically not a tree; however it makes sense for a projected Walsh map to be a tree).

\begin{definition}[Optimal pairing] \label{def:OptimalPairing}
Let $\cB$ be a bubble. A covering of $\cB$ is a $(D+1)$-colored graph obtained from $\cB$ by connecting every white vertex to a black vertex with an edge of color 0. A pairing $\Omega$ is equivalent to a covering denoted $\cB^{\Omega}$ where an edge of color 0 is added between the vertices of each pair of $\Omega$. An optimal pairing is a pairing whose covering maximizes the number of faces of $\cB^{\Omega}$.
\end{definition}

Recall that a projected Walsh map shrinks the reduced submaps $\tilde{\cM}(\cB, \Omega)$ to white vertices. Therefore the covering graphs, which have a single bubble, are mapped to projected maps with a single white vertex. Let $\Omega$ be an optimal pairing. The covering $\cB^{\Omega}$ maximizes the number of faces for graphs with a single bubble. It is mapped to a map $\cT(\cB, \Omega)$ obtained from $\tilde{\cM}(\cB, \Omega)$ by replacing each cilium with an edge carrying the colors incident to its vertex and ending on a univalent black vertex. Therefore, it projects onto a tree: by construction, performing the bijection with an optimal pairing maps the corresponding optimal covering to the unique tree with a single white vertex.

\subsection{Trees}

\begin{figure}
\includegraphics[scale=.4]{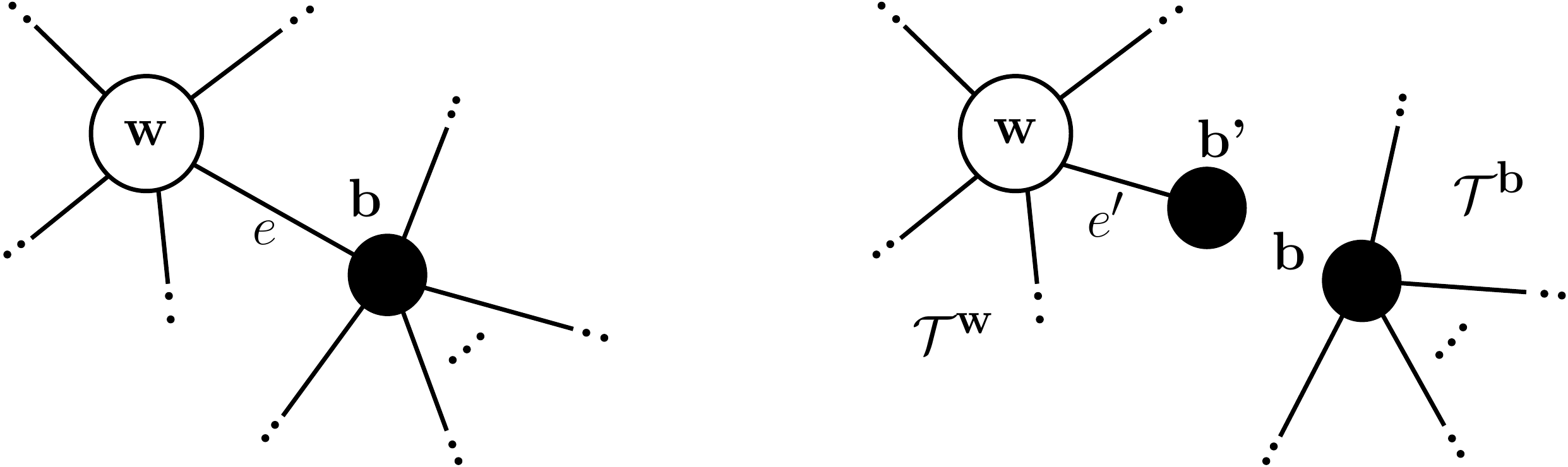}
\caption{\label{fig:Unhooking} A new black vertex is added, and the edge $e$ is unhooked from its black endpoint and attached to the new black vertex.}
\end{figure}

To study projected maps with more than a single white vertex, we will use the following lemma.

\begin{lemma} \label{lemma:Unhooking}
Let $\cW$ be a projected Walsh map and $e$ an edge. We define the \emph{edge unhooking} $\cW\mapsto \cW_e$ of $e$ as in Figure \ref{fig:Unhooking} where a new black vertex is created in the vicinity of a corner incident to $e$ and the edge $e$ is unhooked from its black endpoint and connected to the new black vertex instead ($\cW_e$ may be disconnected).

Denote $\cI(e)$ the set of colors labeling $e$, and $\cI_2(e)\subset \cI(e)$ the set of colors for which two distinct faces run along $e$. Then the variation of the number of faces $\Delta \cF = \cF(\cW_e) - \cF(\cW)$ is
\begin{equation}
\Delta \cF = D-2\lvert\cI_2(e)\rvert.
\end{equation}
\end{lemma}

Notice that if the monochromatic submap $\cW^{(i)}$ is planar, then $i\in\cI_2(e)$ if and only if $e$ is not a bridge in $\cM^{(i)}$.

{\bf Proof.} For each color, there are either 0, 1, or 2 faces running along $e$. The set of colors for which there is no face along $e$ is $\cI_0(e) = [D]\setminus\cI(e)$. Denote $\cI_1(e)\subset \cI(e)$ the set of colors for which there is a single face running along $e$. Clearly $\cI_1(e)$ and $\cI_2(e)$ form a partition of $\cI(e)$ and 
\begin{equation}
|\cI_1(e)| + |\cI_2(e)| = |\cI(e)|.
\end{equation}
For each color $i\in\cI_0(e)$, $\cW_e$ has one more face than $\cW$, sitting around the new black vertex. The face of color $i\in\cI_1(e)$ along $e$ in $\cW$ is split into two faces in $\cW_e$. Conversely the two faces of color $i\in \cI_2(e)$ in $\cW$ are merged in $\cW_e$. This exhausts all the possible colors. Therefore,
\begin{equation}
\Delta\cF = |\cI_0(e)| + |\cI_1(e)| - |\cI_2(e)| = D - 2|\cI_2(e)|.
\end{equation}
\qed

\begin{prop} \label{prop:Trees}
The number of faces of a projected tree $\cT_{\cV_\circ}$ on $\cV_\circ$ white vertices is
\begin{equation}
\cF(\cT_{\cV_\circ}) = (\cF(\cB^\Omega) - D)\cV_\circ + D.
\end{equation}
\end{prop}

{\bf Proof.} One of the white vertices of $\cT_{\cV_\circ}$ is adjacent to a single white vertex, connected to it by a bridge. Unhooking the bridge breaks up the tree into two trees, $\cT_1$ which has a single vertex and $\cT_{\cV_\circ-1}$ with $\cV_\circ - 1$ white vertices. Lemma \ref{lemma:Unhooking} applies and gives
\begin{equation}
\Delta\cF = \cF(\cT_1) + \cF(\cT_{\cV_\circ-1}) - \cF(\cT_{\cV_\circ}) = D,
\end{equation}
since $\cI_2(e) = 0$ for all edges of the tree. A trivial recursion leads to
\begin{equation}
\cF(\cT_{\cV_\circ}) = (\cF(\cT_1) - D)\cV_\circ + D
\end{equation}
and we conclude by noticing that the tree on one white vertex is the covering $\cB^{\Omega}$.
\qed

In the quartic melonic case, $\cF(\cB^\Omega) = 2D-1$ and with $p=2$, we recover the expected relation \eqref{Faces} between the numbers of faces and vertices, $\cF(\cT_{\cV_\circ}) -(D-1)\cV_\circ = D$. Moreover, in that case, $\cF(\cW) \leq D+ (D-1)\cV_\circ$ for a projected map. In a generic model however, it is unclear that adding a new white vertex cannot increase the number of faces more than $(\cF(\cB^\Omega) - D)$, and hence whether $D + (\cF(\cB^\Omega)-D)\cV_\circ - \cF(\cW)$ is always positive or zero. We can nevertheless easily rule out the case where $\cW$ has a single cycle.

\begin{prop}
Let $\Omega$ be an optimal pairing and $\cW$ a projected map with a single cycle. Then $\cW$ has fewer faces than any tree with the same number of white vertices as $\cW$.
\end{prop}

{\bf Proof.} It relies on a lemma proved in \cite{Universality} (lemma 6) which states that in an optimal covering of $\cB$, the number of faces incident to two different edges of color 0 cannot be greater than $D/2$. It can be reformulated as the fact that there cannot be more than $D/2$ faces shared by two pairs of $\Omega$. Since those pairs become edges of projected maps through the bijection, it means that the number of faces shared by two edges incident to a white vertex is bounded by $D/2$. 

Let $\cT$ be a spanning tree in $\cW$ and $e$ the only edge of $\cW$ not in $\cT$, with label set $\cI(e)$, and $\cI_2(e)\subset \cI(e)$ the set of colors for which two distinct faces run along $e$. There is a single path $\cP$ in $\cT$ between the endpoints of $e$, with edges $e_1, \dotsc, e_P$ where $e_1$ is incident to the same white vertex as $e$. Denote $\cI(e, e_1)$ the set of colors of the faces shared by $e$ and $e_1$. Then 
\begin{equation}
\cI_2(e) \subset \cI(e, e_1).
\end{equation}
To prove this, we introduce $\tilde{\cW}$ the stuffed Walsh map which projects on $\cW$, i.e. $P(\tilde{\cW}) = \cW$, and $\tilde{\cW}^{(i)}$ its monochromatic submaps for $i\in[D]$. Let $i\in\cI_2(e)$, then $e$ is not a bridge in $\tilde{\cW}^{(i)}$ (else there would be a single face of color $i$ incident to $e$). It means that there is a path in $\cW^{(i)}$ between the two endpoints of $e$ which does not contain $e$. This path thus has to go along $e_1, \dotsc, e_P$. Therefore $i$ belongs to the set of colors of the faces shared by all those edges and in particular $\cI(e, e_1)$. To conclude, we know from Lemma 6 in \cite{Universality} that
\begin{equation}
|\cI(e, e_1)| \leq D/2,
\end{equation}
and thus 
\begin{equation}
\Delta \cF = D - 2|\cI_2(e)| \geq D - 2|\cI(e, e_1)| \geq 0,
\end{equation}
where the first equality is Lemma \ref{lemma:Unhooking}.
\qed

Finally we can relate the number of faces of a stuffed Walsh map to that of a projected tree in terms of topological quantities.

\begin{theorem} \label{prop:CharDomOrder}
Let $\cW\in\bW_0(\tilde{\cM}(\cB, \Omega))$ and $P(\cW)$ its projection and $\cT$ a projected tree with the same number of white vertices as $P(\cW)$. Then,
\begin{equation}
\cF(\cW) - \cF(\cT) = - D L(\cW) + 2 \sum_{i=1}^D l(\cW^{(i)}) - 2 \sum_{i=1}^D g(\cW^{(i)}).
\end{equation}
Here $l(\cW^{(i)})$ and $g(\cW^{(i)})$ denote the circuit rank and the genus of the monochromatic submap of color $i$ while $L(\cW)$ is the circuit rank of the \emph{projected} map $P(\cW)$.
\end{theorem}

{\bf Proof.} Denote $\cV, \cE, \cF, k, g$ the numbers of vertices, edges, faces, connected components and the genus of a map and $l$ the circuit rank of the underlying graph. For any map, the circuit rank is
\begin{equation}
l = \cE - \cV + k.
\end{equation}
so that Euler's relation reads
\begin{equation} \label{Euler}
\cF = 2k - 2g + \cE - \cV = 2l - 2g - \cE + \cV.
\end{equation}
We denote $\cW^{(i)}$ the monochromatic submap of color $i\in[D]$. Using \eqref{Euler} in each monochromatic submap, the number of faces of $\cW$ can be written
\begin{equation}
\cF(\cW) = \sum_{i=1}^D \cF(\cW^{(i)}) = \sum_{i=1}^D 2l(\cW^{(i)}) - 2g(\cW^{(i)}) - \cE(\cW^{(i)}) + \cV(\cW^{(i)}).
\end{equation}
The formula also works for $\cT$, with the simplification $l(\cT^{(i)}) = g(\cT^{(i)}) = 0$. Moreover, we know from Proposition \ref{prop:Trees} that the number of faces of a projected tree only depends on its number of white vertices and not on its shape. Therefore, we can choose $\cT$ to be obtained out of $\cW$ by unhooking the edges which are not in a spanning tree. This shows that the number of edges with a fixed color label are the same in $\cT$ and $\cW$, hence $\cE(\cW^{(i)}) = \cE(\cT^{(i)})$. Furthermore, every vertex of $\cW$ belongs to all monochromatic submaps, may it be as an isolated vertex, meaning that $\cV(\cW^{(i)}) = \cV(\cW)$ and similarly for $\cT$. Therefore,
\begin{equation}
\cF(\cW) - \cF(\cT) = D \bigl(\cV(\cW) - \cV(\cT)\bigr) + \sum_{i=1}^D 2l(\cW^{(i)}) - 2g(\cW^{(i)}).
\end{equation}
Clearly $\cW$ and $\cT$ have the same number of blue vertices (because they have the same number of submaps $\tilde{\cM}(\cB, \Omega)$), implying that the difference $\cV(\cW) - \cV(\cT)$ is only due to black vertices. This difference is thus preserved by projection of the submaps $\tilde{\cM}(\cB, \Omega)$ onto white vertices, leading to $\cV(\cW) - \cV(\cT) = \cV(P(\cW)) - \cV(P(\cT))$. The projected maps $P(\cW)$ and $P(\cT)$ also have the same number of edges and are both connected, hence
\begin{equation}
\cV(\cW) - \cV(\cT) = \cV(P(\cW)) - \cV(P(\cT)) = l(P(\cT)) - l(P(\cW)) = - l(P(\cW)),
\end{equation}
which concludes the proof.
\qed

\section{Applications} \label{sec:Applications}


\subsection{Melonic graphs}

\begin{figure}
\includegraphics[scale=.7]{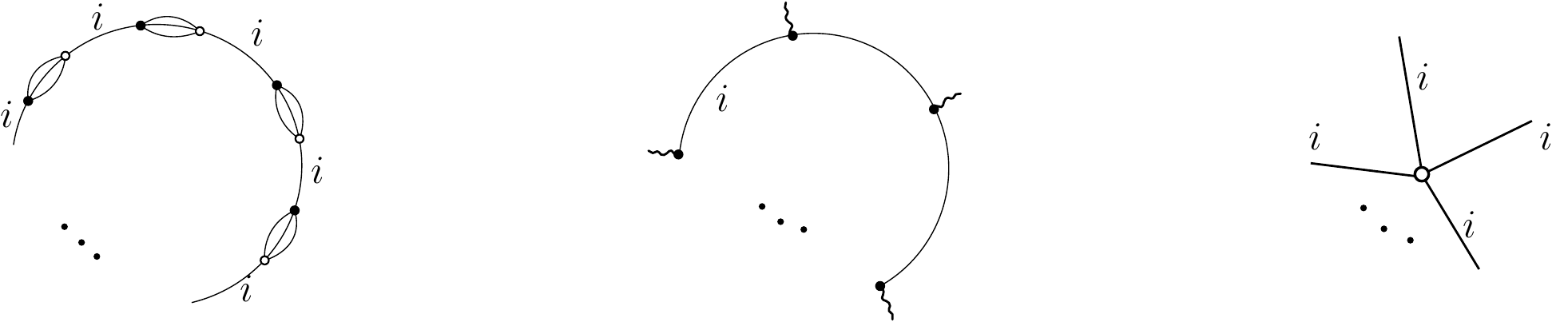} 
\caption{\label{fig:MonoColMel} On the left is a melonic bubble alternating color $i$ and the $D-1$ other colors; next to it, the map $\cM(\cB, \Omega)$ and on the right a reduced map $\tilde{\cM}(\cB, \Omega)$, obtained by replacing the cycle with a star and with cilia replaced with half-edges. Note that the latter is unchanged by the projection onto a white vertex, meaning that in this example stuffed Walsh maps coincide with projected Walsh maps. In fact, it really is a white vertex in the sense of Walsh, i.e. a hyperedge of a hypermap.}
\end{figure}

Consider a \emph{melonic} bubble $\cB$. It can be constructed from the unique bubble on two vertices (with $D$ parallel edges) by adding recursively on any edge a new pair of vertices with $D-1$ parallel edges (an example is shown in Figure \ref{fig:MonoColMel}). This process is called a $(D-1)$-dipole insertion. Vertices thus have a natural pairing $\Omega$ on a melonic bubble, where a pair corresponds to the two vertices added at any dipole insertion. Moreover, $\cB^\Omega$ is the only covering which maximizes the number of faces, making $\Omega$  the only optimal pairing for a melonic bubble \cite{Universality}.

It is easy to see that the edge-colored map $\tilde{\cM}(\cB, \Omega)$ is a tree. Therefore, a stuffed Walsh map is a map whose edges between blue vertices form prescribed subtrees. Notice that with the bubble shown in Figure \ref{fig:MonoColMel}, the stuffed Walsh map really are standard Walsh maps, since then all edges carry the same color and the white vertices can be mapped to hyperedges of fixed degree \cite{Walsh}.

The dominant graphs in $\bG_0(\cB)$ are the melonic ones built out of $\cB$ and they are mapped to projected trees.

\subsection{Bipartite maps}

\begin{figure}
\includegraphics[scale=.7]{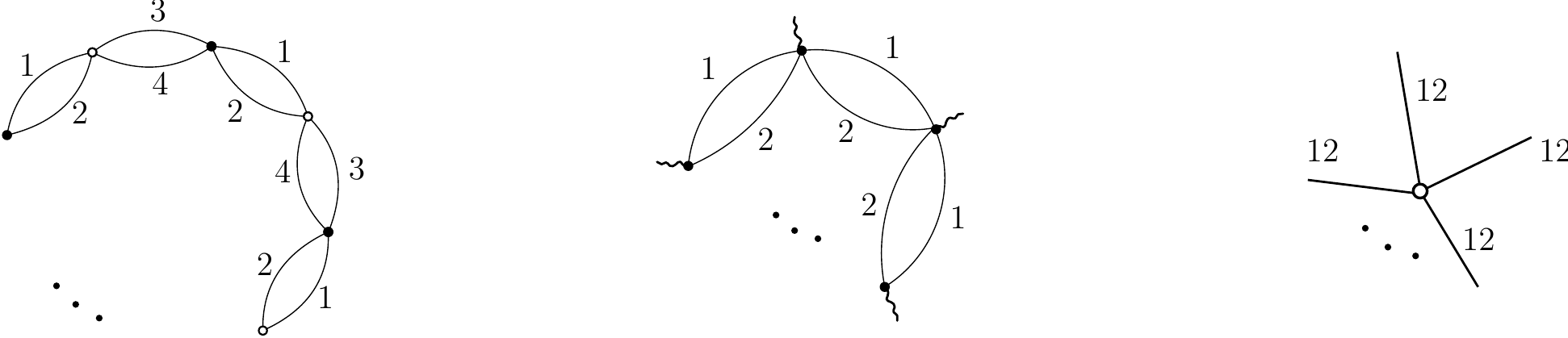} 
\caption{\label{fig:Necklace} On the left, we have a necklace alternating colors 1,2 and colors 3,4. In the middle, this is the map $\cM(\cB, \Omega)$, and on the right a reduced map $\tilde{\cM}(\cB, \Omega)$ with cilia replaced with half-edges, which turns out to coincide with the white vertex of the projected maps.}
\end{figure}

A necklace bubble with $2n$ vertices at $D=4$ is depicted in Figure \ref{fig:Necklace} (there are two others, obtained by transposing the colors 2 and 3, and 2 and 4). All planar coverings are optimal. We choose the pairing which pairs adjacent vertices sharing the edges of colors 3 and 4. The map $\cM(\cB, \Omega)$ thus only has the colors 1 and 2 and is in fact a cycle of parallel edges with colors 1, 2, and cilia on the outside of the cycle. Clearly, we can merge parallel edges into regular edges with the pair of colors $(1, 2)$. To reduce the map so that the monochromatic submaps are forests, we add a vertex inside the cycle and form a star. This reduced map $\tilde{\cM}(\cB, \Omega)$ is shown on the right of Figure \ref{fig:Necklace}. Since $\tilde{\cM}(\cB, \Omega)$ is star-shaped, the stuffed maps coincide with the projected maps. Furthermore, since all edges carry the same couple of colors $(1,2)$, $\tilde{\cM}(\cB, \Omega)$ really is a hyperedge in the sense of Walsh and the stuffed maps are in fact ordinary Walsh maps \cite{Walsh}.

Since each edge is labeled by $D/2=2$ colors, we get from Lemma \ref{lemma:Unhooking} that trees are dominant maps. Dominant maps are thus maps for which the difference $\Delta\cF$ in Theorem \ref{prop:CharDomOrder} vanishes.

For all $\cW\in\bW_0(\tilde{\cM}(\cB, \Omega))$, the monochromatic submaps $\cW^{(1)}=\cW^{(2)}$ coincide with the map $\cW$ itself. This together with the fact that the stuffed maps coincide with their projections, implies that
\begin{equation}
4L(\cW) - 2 \sum_{i=1}^4 l(\cW^{(i)}) = 4l(\cW) - 2 \sum_{i=1}^4 l(\cW^{(i)}) = 0,
\end{equation}
(recall that $L(\cW)$ is the circuit rank of the projected map while $l(\cW)$ is the rank of the map itself). Therefore, Theorem \ref{prop:CharDomOrder} reduces the criterion for dominant maps $\Delta \cF = 0$ to $2g(\cW)=0$, meaning that the dominant maps are the planar maps.

If the three necklace bubbles are now allowed in the colored graphs, each with $2n$ vertices, we choose the pairings so that two vertices of each bubble are paired if they are adjacent and not connected by an edge of color 1. Due to this, all the edges of the maps will have the color 1 in their labels. Again, the stuffed maps and the projected maps are the same. However, they are not ordinary Walsh maps due to the different possible colorings. There are indeed three types of vertices like the one on the right of Figure \ref{fig:Necklace}. The three edges incident to a white vertex all have label $(1, 2)$, or $(1, 3)$, or $(1, 4)$. Notice that $\cW^{(1)} = \cW$.

Trees are dominant maps. Imposing $\Delta\cF=0$ in Theorem \ref{prop:CharDomOrder} then gives 
\begin{equation}
g(\cW^{(i)})=g(\cW)=0 \qquad \text{and} \qquad l(\cW)=l(\cW^{(2)})+l(\cW^{(3)})+l(\cW^{(4)}).
\end{equation}
Consider three forests $\cT^{(i)}$, $i=2, 3, 4$ spanning each $\cW^{(i)}$. The edges in $\cW^{(2)}$, $\cW^{(3)}$, $\cW^{(4)}$ are all distinct so that the edges in the complements $\cW^{(i)}\setminus \cT^{(i)}$ identify $l(\cW^{(2)})+l(\cW^{(3)})+l(\cW^{(4)})$ cycles of $\cW$. Those cycles each have one edge obviously not contained in $\bigcup_{i=2}^4 \cT^{(i)}$. The existence of another cycle, necessarily in $\bigcup_{i=2}^4 \cT^{(i)}$, would then imply $l(\cW)>\sum_{i=2}^4 l(\cW^{(i)})$. The dominant maps are therefore planar maps such that the union $\bigcup_{i=2}^4\cT^{(i)}$ of \emph{any} three forests respectively covering $\cW^{(2)}$, $\cW^{(3)}$ and $\cW^{(4)}$ is a tree. This extends the result of \cite{Melonoplanar} to non-quartic necklace bubbles.

\subsection{A six-vertex graph at $D=4$} \label{sec:6VertexD=4}

\begin{figure}
\includegraphics[scale=.7]{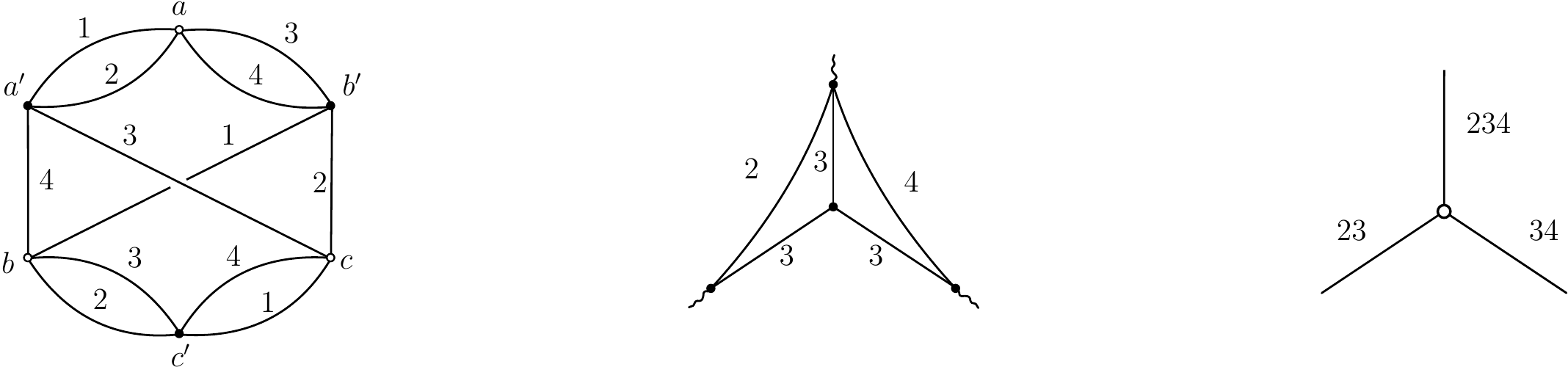} 
\caption{\label{fig:Meander} 
The bubble is on the left. In the center, we have a reduced map for the optimal pairing $\{(a,a'), (b,b'), (c,c')\}$. It can be further simplified to the star on the right which has exactly the same broken faces.}
\end{figure}

We consider the bubble in Figure \ref{fig:Meander}. It has four different optimal pairings, among which $\Omega = \{(a,a'),(b,b'),(c,c')\}$. The map $\cM(\cB, \Omega)$ has a cycle of color 2 and length 2, same for the color 4, and a cycle of length 3 on the color 3. A reduced map is shown in the center of Figure \ref{fig:Meander}, where the cycle of color 3 has been traded for a star, while we have removed an edge of color 2 and one of color 4 too. Looking at the pattern of incidence of broken faces onto the vertices, it is clear that it can be further reduced to a star with different color sets on its three edges as shown on the right of Figure \ref{fig:Meander}. This way, stuffed maps and projected maps are the same.

The white vertex has three distinct edges. Since one of them has the color set $\cI = \{2, 3, 4\}$, hence $|\cI|=3>D/2$, we cannot use directly Lemma \ref{lemma:Unhooking} and conclude that trees are dominant maps. However, the white vertex is such that in a map $\cW$, the submap with edges only labeled by $\cI$ is a forest. It is always possible to complete this forest into a tree $\cT$ spanning $\cW$. The edges of $\cW\setminus \cT$ are all labeled by $D/2=2$ colors. Unhooking them and using Lemma \ref{lemma:Unhooking} shows that
\begin{equation}
\cF(\cW) \leq \cF(\cT).
\end{equation}
As a consequence, trees are dominant maps, and a map $\cW$ is dominant if and only if $\Delta\cF=0$ in Theorem \ref{prop:CharDomOrder}. In this case, since $\cW^{(3)}=\cW$, $\Delta\cF =0$ becomes
\be
\label{eqref:CharDomOrdMean}
l(\cW)-l(\cW^{(2)})-l(\cW^{(4)}) + g(\cW) + g(\cW^{(2)}) + g(\cW^{(4)}) = 0.
\ee
We first analyze the sign of the quantity $l(\cW)-l(\cW^{(2)})-l(\cW^{(4)})$. Consider the forest $\cT^{(234)}$ made of all the edges labeled by $\cI = \{2, 3, 4\}$. It is a sub-forest of $\cW^{(2)}$ and $\cW^{(4)}$. Respectively complete it into two forests, $\cT^{(2)}$ spanning $\cW^{(2)}$ and $\cT^{(4)}$ spanning $\cW^{(4)}$. Since
\begin{equation}
\cW^{(2)} \cup \cW^{(4)} = \cW,
\end{equation}
the union $\cT^{(2)}\cup\cT^{(4)}$ is a connected\footnote{Assume $\cT^{(2)}\cup \cT^{(4)}$ is the disjoint union of two components. Since $\cW$ is connected, there is a path in $\cW$ between them. This path only contains black vertices since a white vertex always has an incident edge of label $\{2, 3, 4\}$ and is thus contained in $\cT^{(2)}\cup \cT^{(4)}$. The path thus has at most a single black vertex $v$, with at least two edges $e_1, e_2$ which connect it to the two components of $\cT^{(2)}\cup \cT^{(4)}$. Without loss of generality, assume that $e_1$ is of type $\{2, 3\}$. It thus belongs to a connected component of $\cW^{(2)}$. The restriction of $\cT^{(2)}$ to this connected component is a spanning tree. Therefore there exists a path in this spanning tree which joins the two ends of $e_1$, meaning that $v$ is in fact connected to one of the two components of $\cT^{(2)}\cup\cT^{(4)}$. We are left with the edge $e_2$ connecting those two components, and a similar reasoning shows that they are in fact connected.} map which spans $\cW$. Let $\cT$ be a spanning tree of $\cT^{(2)}\cup \cT^{(4)}$, hence also a spanning tree of $\cW$. Let $e\in \cW\setminus \cT$. Either $e\in \cT^{(2)}\cup\cT^{(4)}$ or $e\not\in \cT^{(2)}\cup\cT^{(4)}$. In the first case, $e$ lies in $\cT^{(2)}\cup \cT^{(4)}$ but not in its spanning tree and therefore identifies a cycle of $\cT^{(2)}\cup\cT^{(4)}$. If $e\not\in\cT^{(2)}\cup\cT^{(4)}$, it identifies a cycle in $\cW^{(2)}$ or $\cW^{(4)}$. Since the edges carrying both the colors 2 and 4 are labeled with $\{2,3,4\}$ and all belong to $\cT^{(2)}\cup\cT^{(4)}$, it cannot be a cycle common to both $\cW^{(2)}$ and $\cW^{(4)}$. We thus have shown that
\be
l(\cW)=l(\cT^{(2)}\cup\cT^{(4)})+l(\cW^{(2)})+l(\cW^{(4)}).
\ee
In particular, $l(\cW)\geq l(\cW^{(2)})+l(\cW^{(4)})$, so that equation \eqref{eqref:CharDomOrdMean} is equivalent to 
\be
\label{eqref:CharDomOrdMean2}
\left\{
\begin{aligned}
&g(\cW)=g(\cW^{(2)})=g(\cW^{(4)})=0\\
&l(\cW)= l(\cW^{(2)})+l(\cW^{(4)})\quad \Leftrightarrow \quad \text{$\cT^{(2)}\cup\cT^{(4)}$ is a tree.}
\end{aligned}
\right.
\ee


\subsection{The complete bipartite graph in $D=3$} \label{sec:K33}

\begin{figure}
\includegraphics[scale=.7]{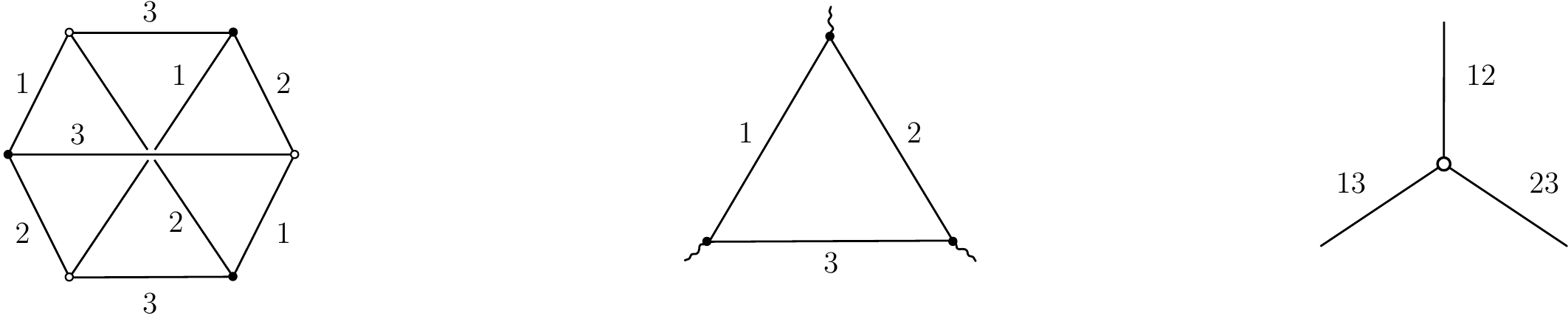} 
\caption{\label{fig:CompBip} The bubble we consider is the complete bipartite graph $K_{3,3}$. In the middle we have drawn a reduced map with $\cB$ as boundary graph. It can be further simplified as shown on the right picture.}
\end{figure}

We consider the complete bipartite graph $K_{3,3}$ with a proper 3-edge-coloring, represented in Figure \ref{fig:CompBip}. It has six different pairings: either vertices inside pairs are all linked by the same color $i=1, 2, 3$, or each one of the three pairs contains a different color. The last three are optimal and equivalent by symmetry. For this pairing $\Omega$, we arrive at a reduced map $\tilde{\cM}(\cB, \Omega)$ like in the middle of Figure \ref{fig:CompBip}. It can be turned into a star-shaped structure without changing the boundary graph, as shown on the right of Figure \ref{fig:CompBip}. This way, the stuffed maps and the projected maps are the same. 

Here, Lemma \ref{lemma:Unhooking} does not directly show that trees are dominant maps. We thus focus on Theorem \ref{prop:CharDomOrder} instead. It gives
\be
\delta(\cT)-\delta(\cW)=3l(\cW)-2\sum_{i=1}^D l(\cW^{(i)})+2\sum_{i=1}^D g(\cW^{(i)}).
\ee
Notice that in contrast with the previous example in Section \ref{sec:6VertexD=4}, no monochromatic submaps $\cW^{(i)}$ cover the whole map $\cW$, because there is no color shared by the three types of edges\footnote{In fact, our approach in Section \ref{sec:6VertexD=4} did not rely on the fact that $\cW^{(3)} = \cW$. It could have been used to first show that $g(\cW) = 0$, implying $g(\cW^{(2)}) = g(\cW^{(4)}) = 0$, leaving us with the analysis of the cycles, which we therefore had to perform anyway.}. However, the union of any two of them covers $\cW$, i.e. for any $i\neq j \in\{1,2,3\}$,
\begin{equation}
\cW^{(i)}\cup\cW^{(j)}=\cW.
\end{equation}
We can therefore perform an analysis of the cycles very similar to that of Section \ref{sec:6VertexD=4}. Indeed, the union of two forests respectively spanning $\cW^{(i)}$ and $\cW^{(j)}$ is a map which spans $\cW$. We choose two such forests in the following way: first observe that $\cW^{(i)}\cap\cW^{(j)}$ is a forest (a cycle would have both colors $(i, j)$ all along which is impossible), which can thus be completed into both a forest $\cT^{(i)}$ spanning $\cW^{(i)}$ and a forest $\cT^{(j)}$ spanning $\cW^{(j)}$. This ensures that $\cW^{(i)}\setminus\cT^{(i)}$ and $\cW^{(j)}\setminus\cT^{(j)}$ do not have any common edges. This choice of spanning forests gives us the following inequalities,
\be
\label{LiLj}
\forall i\neq j \qquad l(\cW)=l(\cT^{(i)}\cup\cT^{(j)})+l(\cW^{(i)})+l(\cW^{(j)})\geq l(\cW^{(i)})+l(\cW^{(j)}).
\ee

Summing the three different relations for $i\neq j\in\{1, 2, 3\}$, we obtain that
\be
3l(\cW)-2\sum_{i=1}^3 l(\cW^{(i)})\geq0. 
\ee
This implies that
\begin{equation}
\Delta \cF = -\Bigl(3l(\cW) - 2 \sum_{i=1}^3 l(\cW^{(i)}) \Bigr) - 2\sum_{i=1}^3 g(\cW^{(i)})
\end{equation}
in Theorem \ref{prop:CharDomOrder} is a sum of integers which are negative or zero. Identifying the maps which maximize the number of faces requires maximizing that quantity. We thus look to find if it is possible to impose all at the same time
\begin{subequations}
\label{K33DomOrd}
\begin{align}
&3l(\cW)=2\sum_{i=1}^3 l(\cW^{(i)}) \label{K33CycleSum}\\
&g(\cW^{(1)})=g(\cW^{(2)})=g(\cW^{(3)})=0 \label{K33Planar}.
\end{align}
\end{subequations}
This is certainly true for trees which consequently are dominant maps. Now we look for the full set of solutions.

Because of the relations \eqref{LiLj}, the constraint \eqref{K33CycleSum} is equivalent to $l(\cW^{(i)}\cup \cW^{(j)})=0$ for all $i\neq j$ which back into \eqref{LiLj} leads to the system
\be
\label{K33VarCond}
l(\cW^{(1)}) + l(\cW^{(2)}) = l(\cW), \qquad l(\cW^{(1)}) + l(\cW^{(3)}) = l(\cW), \qquad l(\cW^{(2)}) + l(\cW^{(3)}) = l(\cW),
\ee
whose solution is
\be
\label{K33VarVarCond}
2l(\cW^{(1)}) = 2l(\cW^{(2)}) = 2l(\cW^{(3)}) = l(\cW).
\ee

\begin{figure}
\includegraphics[scale=1.3]{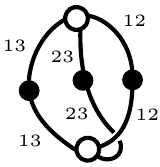}\hspace{3cm} \includegraphics[scale=1.5]{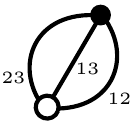} \hspace{3cm}\includegraphics[scale=1.5]{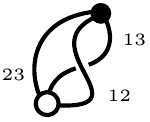} 
\caption{\label{fig:DomCells} This figure shows the bridgeless dominant maps with two fundamental cycles.}
\end{figure}

We now study the case $l(\cW)=2$, for which \eqref{K33VarVarCond} rewrites
\be
l(\cW^{(1)})=l(\cW^{(2)})=l(\cW^{(3)})=1,
\ee
and identify the bridgeless solutions. There are three distinct cycles $C_1, C_2, C_3$, one for each color, but only two fundamental cycles. It means that every cycle is the symmetric difference of the other two. $\cW$ thus has the structure of a Theta graph, i.e. two nodes with three segments between them. Each segment must be part of two cycles, meaning that the edges of a segment all have the same couple of colors. Bipartiteness and the structure of the white vertex prevent any map from having a chain with more than two consecutive edges with the same couple of colors. Therefore, each segment has one or two edges, with the same colors. The allowed maps are thus restricted to those shown in Figure \ref{fig:DomCells}.

We now prove by induction on the circuit rank $l$ that any map $\cW$ which has no submap homeomorphic to one of those of Figure \ref{fig:DomCells} and which is not a tree verifies
\be
\label{eqref:InequK33}
3l(\cW)-2\sum_{i=1}^3 l(\cW^{(i)})>0,
\ee
and is therefore not a dominant map.

We saw that this property is true for $l(\cW)\leq 2$. Let $l>2$ and $\cW$ such that $l(\cW)=l$. We distinguish two cases.

First, assume there exists an edge $e$ which is not a bridge and such that $\lvert\cI_2(e)\rvert\leq 1$. By Lemma \ref{lemma:Unhooking}, unhooking this edge gives rise to a map with more faces, so $\cW$ cannot be dominant.

If there is no such edge, then all edges which are not bridges are such that $|\cI_2(e)| = 2$. Let $e_{ij}$ be one of them, with colors $\{i,j\}$. It is not a bridge in any of the two monochromatic submaps it is contained in. Let us unhook that edge, which leads to a map $\cW'$ with $l(\cW') = l(\cW) - 1$ and $\sum_{i=1}^3 l(\cW'^{(i)}) = \sum_{i=1}^3 l(\cW^{(i)}) - 2$, so that $3 l(\cW) - 2 \sum_{i=1}^3 l(\cW^{(i)}) = 3 l(\cW') - 2\sum_{i=1}^3 l(\cW'^{(i)}) - 1$.

As $e_{ij}$ is not a bridge in $\cW^{(i)}$ and $\cW^{(j)}$, the adjacent edges $e_{jk}, e_{ik}$ at its white endpoint also are not bridges respectively in $\cW^{(j)}$ and $\cW^{(i)}$. Indeed, $e_{ij}$ belongs to two distinct cycles, one containing $e_{ik}$ and whose edges all contain the color $i$ and one containing $e_{jk}$ and whose edges all contain the color $j$. Let us denote those two cycles $C_i=(e_{ik}, \dotsc, e_{ij})$, and $C_j=(e_{jk}, \dotsc, e_{ij})$, and $e'_{ij}$ the first edge they have in common (which might be $e_{ij}$). Then the concatenation of the two chains $(e_{ik}, ...,e'_{ij})\subset C_i $ and $(e_{jk}, ...e'_{ij})\subset C_j $, with $e'_{ij}$ excluded, is also a cycle in $\cW'$. This implies that $e_{ik}$ and $e_{jk}$ are not bridges in $\cW'$. Furthermore, $e_{ij}$ is now attached to a leaf in $\cW'$, so that $e_{ik}$ and $e_{jk}$ are respectively bridges in $\cW'^{(i)}$ and $\cW'^{(j)}$. When unhooking, say, $e_{ik}$ from its black endpoint, we get a map $\cW''$ such that $l(\cW'')=l(\cW')-1$ and $\sum_{i=1}^3 l(\cW''^{(i)})=\sum_{i=1}^3 l(\cW'^{(i)})-\lambda$, with $\lambda\in\{0,1\}$. It comes that
\begin{equation}
3l(\cW)-\sum_{i=1}^3 l(\cW^{(i)}) = 3 l(\cW') - 2\sum_{i=1}^3 l(\cW'^{(i)}) - 1 = 3l(\cW'')-2\sum_{i=1}^3 l(\cW''^{(i)})+2(1-\lambda) \geq 3l(\cW'')-2\sum_{i=1}^3 l(\cW''^{(i)}).
\end{equation}
Notice that $\cW''$ is not a tree as $l(\cW)>2 \Rightarrow l(\cW'')>0$. Moreover, unhooking $e_{ij}$ and $e_{ik}$ cannot create a submap homeomorphic to one of Figure \ref{fig:DomCells}. The induction hypothesis thus applies on $\cW''$ from which we conclude that $3l(\cW)-\sum_{i=1}^3 l(\cW^{(i)}) >0$.

\begin{figure}
\includegraphics[scale=1.2]{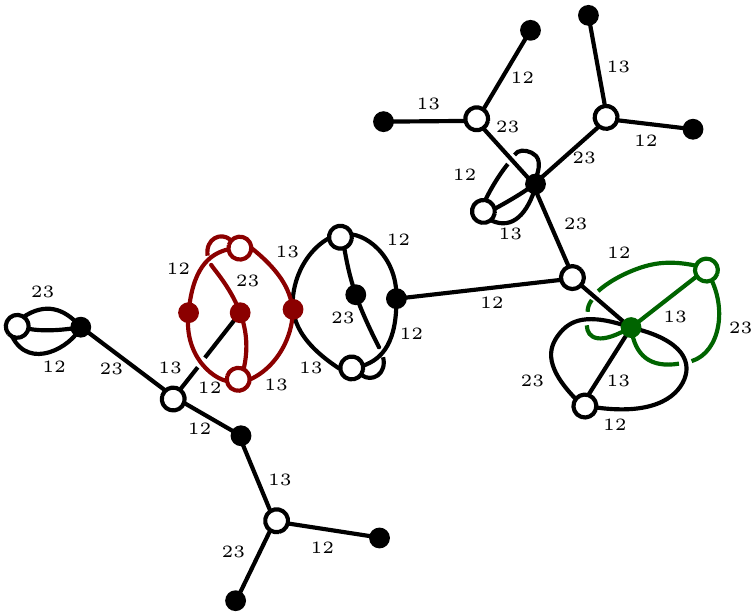} 
\caption{\label{fig:ExDomOrdK33} This is an example of a dominant map.}
\end{figure}

As a consequence, a dominant map $\cW$ has its monochromatic submaps planar and is characterized as follows. First break it into 1-edge-connected components, then identify the submaps homeomorphic to Figure \ref{fig:DomCells} and unhook two of their edges, and then repeat this process until $\cW$ is a tree. An example is given in Figure \ref{fig:ExDomOrdK33}.

\section{Matrix models for stuffed Walsh maps} \label{sec:MatrixModel}

$(D+1)$-edge-colored graphs with bubbles in a finite set $\{\cB_\alpha\}_{\alpha\in A}$ are generated by a \emph{tensor model}. First assign a polynomial invariant under $U(N)^D$ to each bubble in the following way. Every white vertex of $\cB$ receives a copy of a tensor ${\bf T}$ and every black vertex its complex conjugate ${\bf \overline{T}}$. This tensor has $D$ indices, each ranging from 1 to $N$, and components $T_{a_1 \dotsb a_D}$. When an edge of color $i$ connects two vertices in $\cB$, the indices in position $i$ of their corresponding tensors are identified and summed over, $\sum_{a_i, a'_i=1}^D T_{a_1\dotsb a_i\dotsb a_D} \bar{T}_{a'_1 \dotsb a'_i \dotsb a'_D} \delta_{a_i, a'_i}$. Performing all those contractions for a bubble, one obtains an invariant polynomial denoted $\cB({\bf T}, {\bf \overline{T}})$. In particular, there is a single bubble on two vertices, which corresponds to the single quadratic invariant ${\bf T}\cdot \overline{{\bf T}} = \sum_{a_1, \dotsc, a_D} T_{a_1\dotsb a_D} \overline{T}_{a_1\dotsb a_D}$.

Consider
\begin{equation}
\label{eqref:ArbTensZ}
Z(\{\lambda_\alpha\},N) = \exp\,F = \int \exp\left( -\sum_{\alpha=1}^A N^{s_{\alpha}} \lambda_\alpha \cB_\alpha({\bf T}, \overline{{\bf T}})\right)\, d\mu_0({\bf T}, \overline{{\bf T}}),
\end{equation}
where $d\mu_0({\bf T}, \overline{{\bf T}})$ is the Gaussian measure,
\be
\label{eqref:gaussianmeasure}
d\mu_0({\bf T}, \overline{{\bf T}}) = \exp \left(- N^{D-1} {\bf T}\cdot \overline{{\bf T}}\right)\ \prod_{a_1, \dotsc, a_D} \frac{dT_{a_1 \dotsb a_D} d\overline{T}_{a_1 \dotsb a_D}}{2i \pi}.
\ee
We perform a Feynman expansion of $Z$. First expand each $\exp(N^{s_{\alpha}} \lambda_\alpha \cB_\alpha({\bf T}, \overline{{\bf T}}))	$ as a series in $\lambda_\alpha$ and commute the sums with the integral. We are then left with Gaussian moments. By Wick's theorem, they are evaluated as sums over all pairings between $T$'s and $\overline{T}$'s. Those pairings are represented by new edges to which the color 0 is assigned. Thus $F = \ln Z$ generates connected $(D+1)$-colored graphs whose bubbles are from the set $\{\cB_\alpha\}_{\alpha\in A}$.

In the present section, we rewrite tensor models as matrix models which generate stuffed Walsh maps. We will focus on the case  of a single bubble $\cB$ with partition function
\be
Z_{\cB}(\lambda, N)=\int \exp\biggl[-\lambda N^s \cB({\bf T}, {\bf \bar{T}})\biggr]d\mu_{0}({\bf T}, {\bf \bar{T}}).
\ee
The generating function of the cumulants is 
\be
\ln Z_{B}[\lambda, N;{\bf J, \bar J}]
=
\ln\int e^{-{\lambda}{N^s} \cB({\bf T}, {\bf \bar{T}})-{\bf \bar{J}.T-\bar{T}.J}}d\mu_{0}({\bf T}, {\bf \bar{T}}).
\ee
The matrix model we are going to present is a \emph{multi-matrix model, with matrices of various sizes, and multi-trace interactions}. It has a Gaussian measure, while its potential is a sum of two distinct terms. The first one stands for the interaction bubble $\cB$ and can be graphically interpreted as $\tilde{\cM}(\cB, \Omega)$, while the second one, which does not depend on $\cB$, is an infinite series in the conjugate matrices which performs all the possible gluings of interaction bubbles and corresponds to the black vertices of stuffed Walsh maps.

The matrices involved in this new representation are indexed by all possible color sets $\cI\in\cP([D])$, where $\cP([D])$ is the set of subsets of $[D]$. The matrix $\sigma_{\cI}$ is of size $N^{\lvert\cI\rvert}\times N^{\lvert\cI\rvert}$, and has complex entries $\sigma_{\cI\mid(i_{1},\dotsc, i_{\lvert\cI\rvert});(j_{1},\dotsc, j_{\lvert\cI\rvert})}$.

We introduce the notations $\hI= [D]\setminus\cI$, and if $\cI=\{k(1),\dotsc,k(\lvert\cI\rvert)\}$, then $\fsig_{\cI\mid (i_{1},\dotsc, i_{D});(j_{1},\dotsc, j_{D})}$ will stand for $(\prod_{k\in\hI}\delta_{i_{k},j_{k}})\times\sigma_{\cI\mid (i_{k(1)},\dotsc, i_{k( \lvert\cI\rvert )});(j_{k(1)},\dotsc, j_{k(\lvert\cI\rvert)})}$. Note that while the $\sigma_{\cI}$ generally do not commute, $[\fsig_{\cI},\fsig_{\cJ}]=0$ if $\cI\cap\cJ=\emptyset$.

\subsection{Matrix bubbles}

Given a bubble $\cB$ with $\cV$ black vertices and an arbitrary labeling of both its black and white vertices from $1$ to $\cV$, there exist $D$ permutations $\tau_1, \dotsc, \tau_D$ of $\{1, \dotsc, \cV\}$ which are defined by the fact that $\tau_i$ maps the white vertex of label $a$ to the black vertex of label $\tau_i(a)$ if there is an edge of color $i$ between $a$ and $\tau_i(a)$. Then $\cB$ can be written in the following form
\be
\cB({\bf T}, {\bf \bar T})= 
 \sum_{  \substack {  { i^1_{1},\dotsc, i^1_{D},\ \   }   \\ {\dotsc\ \ \ }  \\ { i^{\cV}_{1},\dotsc, i^{\cV}_{D}=1  }  }   }^N
 \sum_{  \substack {  { j^1_{1},\dotsc,j^1_{D},\ \   }   \\ {\dotsc\ \ \ }  \\ { j^{\cV}_{1},\dotsc, j^{\cV}_{D}=1  }  }   }^N  
 \prod_{a,b=1}^{\cV}\ T^{}_{i^a_1 \dotsb i^a_D}\bar T^{}_{j^b_{1} \dotsb j^b_{D}}
 \prod_{k=1}^D\delta_{i_k^a, j_k^{\tau_{k}(a)}}.
\ee

Since the key ingredient to go from $\cB$ to the edge-colored map $\tilde{\cM}(\cB, \Omega)$ which encodes $\cB$ through its broken faces is a pairing $\Omega$ of $\cB$, we re-organize the sums over tensor indices according to $\Omega$. The pairing can also be thought of as a permutation $\tau_0$ that maps each white vertex $a$ to its black partner $\tau_0(a)$. The permutation $\tau_0$ induces a partition of the edges of $\cB$ as follows. For each white vertex $a$, we partition the $D$ edges incident to $a$ according to whether they are incident to its black partner $\tau_0(a)$ or to another black vertex. Denote $\hI_a$ the set of colors incident to both $a$ and $\tau_0(a)$, i.e. if $i\in\hI_a$ then $\tau_i(a)=\tau_0(a)$. Its complement $\cI_a = [D]\setminus \hI_a$ is the set of colors incident to $a$ but not to $\tau_0(a)$, i.e. for all $i\in\cI_a$, $\tau_i(a)\neq \tau_0(a)$.

We re-organize the above expression of $\cB$ according to that partition,
\be
\cB({\bf T}, {\bf \bar T})=
\sum_{  \substack {  { i^1_{1},\dotsc, i^1_{D},\ \   }   \\ {\dotsc\ \ \ }  \\ { i^{\cV}_{1},\dotsc, i^{\cV}_{D}=1  }  }   }^N\ 
 \sum_{  \substack {  { j^1_{1},\dotsc, j^1_{D},\ \   }   \\ {\dotsc\ \ \ }  \\ { j^{\cV}_{1},\dotsc, j^{\cV}_{D}=1  }  }   }^N
 \prod_{a=1}^{\cV}\ \ \bigl[ T_{i^a_1 \dotsb i^a_D}\bar T_{j^{\tau_0(a)}_{1}\dotsb j^{\tau_0(a)}_{D}}
 \prod_{p\in \hI_a} \delta_{i_p^a, j_p^{\tau_p(a)}} \bigr] 
 \prod_{q\in \cI_a} \delta_{i_{q}^a, j_{q}^{\tau_q(a)}}.
\ee
Each index $i_c^a$ on the color $c$ and incident to $a$ is identified with the index $j^{\tau_c(a)}_c$ incident to the black vertex $\tau_c(a)$. We first perform the sum over those indices whose colors are in the sets $\hI_a$, i.e. which connect the vertices of a pair of $\Omega$. After performing this sum, we recognize the matrix elements of the partial contraction of ${\bf T}$ and ${\bf \bar T}$ along the set of colors $\hI_a$, denoted ${\bf T}\cdot_{\hI_{a}}\bar {\bf T}$. The remaining sum runs over the indices $i_c^a$ of colors $c\in\cI_a$. Denote $\cI_a = \{  q_{a,1}, ... ,q_{a,\lvert\cI_a\rvert}  \}$. It comes 
\be
\label{bubble}
\cB({\bf T}, {\bf \bar T}) =
\sum_ {   
\substack {
{   \{i^\alpha_{q_{\alpha,1}},\dotsc, i^\alpha_{q_{\alpha,\lvert\cI_\alpha\rvert}}\ \   }\\
{   j^\alpha_{q_{\alpha,1}},\dotsc, j^\alpha_{q_{\alpha,\lvert\cI_\alpha\rvert}}\}   }
}
}\ 
 \prod_{a=1}^{\cV}\ \bigl[ {\bf T}\cdot_{\hI_{a}}\bar {\bf T}\bigr]_{ (i^a_{q_{a,1}},\dotsc, i^a_{q_{a,\lvert\cI_a\rvert}}); (j^a_{q_{a,1}},\dotsc, j^a_{q_{a,\lvert\cI_a\rvert}})    } \prod_{q\in\cI_a  } \delta_{i_q^a, j_q^{\tau_q(a)}},
\ee

In order to obtain the matrix potential of the bubble $\cB$, we simply use the above expression with the replacement rule 
\begin{equation}
{\bf T}\cdot_{\hI_a} {\bf \bar T} \rightarrow \sigma_{\cI_a}.
\end{equation}
Each pairing $\Omega$ thus defines a matrix potential for $\cB$,
\be
\label{matrixbubble}
V_{\cB, \Omega}^\circ(\{ \sigma_{\cI}\})=
\sum_{   
\substack {
{   \{i^\alpha_1,\dotsc, i^\alpha_{\lvert\cI_\alpha\rvert} \ \   }\\
{   j^\alpha_1,\dotsc, j^\alpha_{\lvert\cI_\alpha\rvert}\}   }
}}\ 
 \prod_{a=1}^{\cV}\sigma_{\cI_a\mid (i^a_1,\dotsc, i^a_{\lvert\cI_a\rvert}); (j^a_1,\dotsc, j^a_{\lvert\cI_a\rvert})    } \prod_{q\in\cI_a  } \delta_{i_q^a, j_q^{\tau_q(a)}}.
\ee
Clearly, there is one matrix $\sigma_\cI$ for each couple of vertices of $\cB$ paired in $\Omega$ and an index of color $i$ is contracted between two matrices if and only if there is a broken face of color $i$ between the corresponding pairs in $\cM(\cB, \Omega)$. We therefore have an algebraic representation of $\cM(\cB, \Omega)$.

\subsection{The matrix models}

\begin{theorem}
\label{thm:intmat1}
For any bubble $\cB$ equipped with a pairing $\Omega$, the partition function and generating function of cumulants rewrite as
\begin{gather} \label{TensorMatrixEquality}
 Z_{B}(\lambda, N)
=\int_{\bC^D}e^{-{\lambda}{N^s} \cB({\bf T}, {\bf \bar{T}})}d\mu_{0}({\bf T}, {\bf \bar{T}})
=\int e^{ - {\lambda}{N^s} V_{\cB,\Omega}^{\circ}(\{\sigma_{\cI}\})-V^{\bullet}(\{\bar{\sigma}_\cI\})  }\prod_{\cI}d\mu_{0}(\sigma_\cI,\bar{\sigma}_\cI),
\\
\label{eqref:genfunc}
 Z_{B}[\lambda, N;{\bf J, \bar J}]
=\int e^{ - {\lambda}{N^s} V_{\cB,\Omega}^{\circ}(\{\sigma_{\cI}\})-V^{\bullet}(\{\bar{\sigma}_\cI\})-\Tr \bigl[({\bf J}\otimes{\bf \bar{J}}).(\un^{\otimes D}+\sum_{\cI\in\cP([D])} \bar \fsig_\cI)^{-1} \bigl] }\prod_{\cI}d\mu_{0}(\sigma_\cI,\bar{\sigma}_\cI),
\end{gather}
where $V_{\cB,\Omega}^{\circ}(\{\sigma_\cI\})$ is given in \eqref{matrixbubble}, $d\mu_{0}(\sigma_\cI,{\bf  \bar \sigma_\cI})$ is the Gaussian measure which propagates each matrix element of $\sigma_\cI$ to the same element of $\bar{\sigma}_\cI$, and
\be
\label{BlackVert}
V_{\bullet}(\{\bar{\sigma}_\cI\}) =\Tr\log\biggl(\un^{\otimes D}+\sum_{\cI\in\cP([D])}
\bar \fsig_\cI\biggr).
\ee
These equalities between tensor integrals and matrix integrals are \emph{perturbative}, i.e. they hold order by order in powers of the coupling constant $\lambda$.
\end{theorem}

The strategy is simply to expand the exponentials of the non-quadratic terms and prove the equalities order by order in $\lambda$. The proof then relies on the one-dimensional case which we present as a lemma.
\begin{lemma} 
For all $a\in\bC$ and $p\in\bN$, the following relation holds
\be
\label{eqref:relationprf}
\langle z^p\,e^{-a\bar{z}}\rangle_0 \equiv \int_{\mathbb C} z^p \ e^{-  a \bar z }d\mu_{0}(z, \bar z) = a^p.
\ee
\end{lemma}
\prf 
\be 
\int_{\bC} z^p \ e^{-  a\bar z -z\bar z}\frac{dz d\bar z}{\pi}
=\biggl[ e^{\frac{\partial}{\partial z}\frac{\partial}{\partial \bar z}} e^{a\bar z}z^p\biggr]_{z=\bar z=0}
=\sum_{k\ge0}\frac{1}{k!}\bigl[ \bigl(\frac{\partial}{\partial z}\bigr)^k z^p\bigr]_{z=0}\bigl[ \bigl(\frac{\partial}{\partial \bar z}\bigr)^k e^{a\bar z}\bigr]_{\bar z=0}
=\sum_{k\ge0}\frac{1}{k!}(p!\delta_{p,k})(a^k)=a^p . 
\ee
\qed

This is easily extended to products of multivariate polynomials.
\begin{lemma} 
\label{polprod}
For any $K\in\bN,\ \{a_i\}\in\bC^K$, and given  $P[z_1,\dotsc,z_K]$ and $Q[z_1,\dotsc,z_K]$ two complex polynomials, the following relation holds
\begin{equation}
\begin{aligned}
\langle P(z_1,\dotsc,z_K) Q(z_1,\dotsc,z_K)\ e^{- \sum_i a_i\bar z_i } \rangle_0 &\equiv \int_{\mathbb C ^K} P(z_1,\dotsc,z_K) Q(z_1,\dotsc,z_K)\ e^{- \sum_i a_i\bar z_i }\prod_id\mu_{0}(z_i, \bar z_i)\\
&= P(a_1,\dotsc, a_K) Q(a_1,\dotsc, a_K).
\end{aligned}
\end{equation}	
\end{lemma}

\prf The generalization of \eqref{eqref:relationprf} to $K$ complex variables $z_1,..,z_K$ is the obvious relation
\be
\label{eqref:monomial}
\int_{\mathbb C ^K} \prod_{i=1}^K z_i^{p_i} \ e^{- \sum_i a_i.\bar z_i }\prod_id\mu_{0}(z_i, \bar z_i)=\prod_{i=1}^Ka_i^{p_i},
\ee
so that by linearity, for two polynomials $P(z_1,\dotsc,z_K)=\sum_{m}\alpha_m\prod_i z_i^{p_{m,i}}$ and $Q(z_1,\dotsc,z_K)=\sum_{n}\beta_n\prod_i z_i^{q_{n,i}}$, 
\begin{multline}
\langle P(z_1,\dotsc,z_K) Q(z_1,\dotsc,z_K)\ e^{- \sum_i a_i\bar z_i } \rangle_0\\
=\sum_{m,n}\alpha_m\beta_n\langle \prod_i z_i^{p_{m,i}+q_{n,i}}e^{-\vec{a}.\vec{\bar z}}\rangle_0=\sum_{m,n}\alpha_m\beta_n\langle \prod_i a_i^{p_{m,i}+q_{n,i}}e^{-\vec{a}.\vec{\bar z}}\rangle_0=P(a_1,..,a_K)Q(a_1,..,a_K)\nonumber . 
\end{multline}
\qed

{\bf Proof of Theorem \ref{thm:intmat1}.} Given a bubble $\cB$, $\Omega$ a pairing of its vertices, and $q\in\bN$, we now prove that the following equality holds
\be
\label{bubblepower}
\bigl[ \cB({\bf T},{\bf \bar T}) \bigr]^p
= 
\int  \bigl[V_{\cB,\Omega}^{\circ}(\{\sigma_{\cI}\})\bigr]^p\ \exp\biggl[ - \sum_{\cI}  \Tr\ [{\bf T}\cdot_{\hI}{\bf \bar T}]\sigma_\cI^\dagger  \biggr]
\prod_{\cI}d\mu_{0}(\sigma_\cI,{\bf \sigma_\cI^\dagger}).
\ee

We use expression \eqref{bubble} for the bubble $\cB$. Relation \eqref{eqref:monomial} applied to the products over the vertices of matrix elements of ${\bf T}\cdot_{\hI_{a}}\bar {\bf T}$ gives
\begin{multline}
\prod_{a=1}^{\cV}\ [ {\bf T}\cdot_{\hI_{a}}\bar {\bf T}\bigr]_{ (i^a_{q_{a,1}},\dotsc, i^a_{q_{a,\lvert\cI_a\rvert}}); (j^a_{q_{a,1}},\dotsc, j^a_{q_{a,\lvert\cI_a\rvert}}) } \\
=
\int  \prod_{a=1}^{\cV}\ \sigma_{\cI_a\mid (i^a_{q_{a,1}},\dotsc, i^a_{q_{a,\lvert\cI_a\rvert}}); (j^a_{q_{a,1}},\dotsc, j^a_{q_{a,\lvert\cI_a\rvert}}) } 
\ e^{ - \sum_{\cI}  \Tr\bigl( [{\bf T}\cdot_{\
\hI}{\bf \bar T}] \sigma_\cI^\dagger\bigr)  }
\prod_{\cI}d\mu_{0}(\sigma_\cI,{\bf  \sigma_\cI^\dagger}),
\end{multline}
since the integration over the matrix elements of $\sigma_\cI$ for $\cI\neq \cI_a$ contributes to a factor one, as well as integration over the matrix elements of $\sigma_{\cI_a}$ with indices different from $(i^a_{q_{a,1}},\dotsc, i^a_{q_{a,\lvert\cI_a\rvert}}); (j^a_{q_{a,1}},\dotsc, j^a_{q_{a,\lvert\cI_a\rvert}})$.

Inserting this relation into \eqref{bubble} and using expression \eqref{matrixbubble}, one gets relation \eqref{bubblepower} for $p=1$. The relation for $p\ge0$ then follows from applying Lemma $\ref{polprod}$.

Now notice that $\Tr\ [{\bf T}\cdot_{\hI}{\bf \bar T}]\sigma_\cI^\dagger = \sum_{ \substack {   {i_1,\dotsc, i_D} \\ {j_1,\dotsc, j_D}  } }T_{i_1\dotsb i_D}\bar T_{j_1\dotsb j_D}\bar\fsig_{\cI\mid(i_1,\dotsc, i_D);(j_1,\dotsc, j_D)}$ so that 
\be
\int d\mu_{0}({\bf T}, {\bf \bar{T}})\exp\biggl[ - \sum_{\cI}  \Tr [{\bf T}\cdot_{\hI}{\bf \bar T}]\sigma_\cI ^\dagger \biggr]
=\frac{1}{\det\bigl[ \un^{\otimes D}+\sum_{\cI}\bar\fsig_\cI  \bigr]}
=\exp\bigl[ -\Tr \ln (\un^{\otimes D}+\sum_{\cI}\bar\fsig_\cI ) \bigr].
\ee

Integrating \eqref{bubblepower} over ${\bf T,\bar T}$ thus yields the matching of each term of the perturbative series
\be
\label{eqref:pBubbles}
\frac{(-N^s\lambda)^p}{p!}\int \bigl[ \cB({\bf T},{\bf \bar T}) \bigr]^p d\mu_{0}({\bf T}, {\bf \bar{T}})
=
\frac{(-N^s\lambda)^p}{p!}\int  \bigl[V_{\cB,\Omega}^{\circ}(\{\sigma_{\cI_c}\})\bigr]^p\ e^{ -\Tr \ln (\un^{\otimes D}+\sum_{\cI}\bar\fsig_\cI ) }
\prod_{\{\cI\}}d\mu_{0}(\sigma_{\cI},{\bf \bar \sigma_{\cI}}).
\ee
Equation \eqref{eqref:genfunc} works similarly with sources.
\qed

We now perform the Feynman expansion of the matrix model \eqref{TensorMatrixEquality}. There are two types of non-quadratic terms in the action which give rise to two types of vertices. 
\begin{itemize}
\item $V^\bullet$ in \eqref{BlackVert}, is expanded as an infinite sum of monomials,
\be
V^\bullet(\{\bar\sigma_\cI\})
=\Tr\ln(\un^{\otimes D}+\sum_\cI\bar\fsig_\cI)=\Tr\sum_{k>0}\frac{(-1)^{k+1}}{k}(\sum_\cI\bar\fsig_\cI)^k
=\sum_{k>0}\frac{(-1)^{k+1}}{k}\sum_{\cI_1,..,\cI_k}\Tr\bar\fsig_{\cI_1}\dotsm \bar\fsig_{\cI_k},
\ee
where the last sum is taken over non-necessarily distinct $\cI_j$ and up to a cyclic permutation of the sets. To each distinct monomial $\Tr\bar\fsig_{\cI_1}\dotsm \bar\fsig_{\cI_k}$, one associates a \emph{black} vertex of degree $k$. Its incident half-edges are labeled with the color sets $\cI_1,\dotsc, \cI_k$. They are cyclically ordered (clockwise or counter-clockwise) according to their order in the trace.

Since a matrix $\fsig_\cI$ represents a pair of vertices of $\cB$, it is clear that $\Tr\bar\fsig_{\cI_1}\dotsm \bar\fsig_{\cI_k}$ corresponds to a cycle of a colored graph which alternates pairs of vertices and edges of color 0 (represented as corners around the vertex). One thus recovers the black vertex of stuffed Walsh maps, introduced in Section \ref{sec:BlackVertex}.

\item To $V^{\circ}_{\cB,\Omega}$, one can associate a \emph{white} vertex of degree $\cV$ (half the number of vertices of the bubble $\cB$). The incident edges are labeled with color sets $\cI$ which correspond to the labels of the matrices $\sigma_\cI$ which do appear in the expression of $V^\circ_{\cB, \Omega}$. The order of those incident edges around each white vertex is fixed.

\item The quadratic term of the action propagates a $\sigma_\cI$ to a $\bar{\sigma}_{\cI}$ and thus gives rise to edges which carry a color set $\cI\in \cP([D])$ and connect black to white vertices.

\item When $V^{\circ}_{\cB,\Omega}$ is represented as a white vertex, we thus recover the \emph{projected} Walsh maps. However, in matrix models (and tensor models), a map of the Feynman expansion receives a factor $N$ for each face of each color. Clearly, representing $V^{\circ}_{\cB,\Omega}$ does not allow to track the faces graphically. To do so, one needs to represent graphically the pattern of identification of matrix indices in $V^{\circ}_{\cB,\Omega}$. By construction, this is what the edge-colored map $\cM(\cB, \Omega)$ (or $\tilde{\cM}(\cB \Omega)$) does. This way, one sees that the matrix model of Theorem \ref{thm:intmat1} generates the stuffed Walsh maps $\bW_0(\cM(\cB, \Omega))$.
\end{itemize}

\section*{Conclusion}

We have established a bijection between colored graphs with prescribed bubbles and stuffed Walsh maps. Those maps extend the ordinary hypermaps, which can be represented as bipartite maps \cite{Walsh}, by stuffing the hyperedges with prescribed maps. This shows in a precise way to what extent edge-colored graphs generalize combinatorial maps and in particular the recently introduced stuffed maps \cite{Stuffed}.

We hope that this bijection will be useful to tackle some interesting problems about colored graphs. The main objective is to classify the graphs according to their number of faces at fixed number of vertices or bubbles. This was done in \cite{GS} for generic colored graphs (no prescription on the allowed bubbles) and with a different notion of faces (including those which carry couples of colors $\{i,j\}$ for $i,j\neq 0$). If $\cB$ is a melonic bubble with $p$ white vertices, then $\bG_0(\cB)$ contains melonic graphs for which we know the maximal ratio $r = (D-1)(p-1)$ between the number of faces and the number of bubbles. It is then possible to perform an analysis similar to that \cite{GS}, see \cite{DSSD}. However, if $\cB$ is generic, the maximal value of the ratio $r$ between the number of faces and the number of bubbles is unknown. Using the bijection, we have found in Section \ref{sec:DominantMaps} that projected trees have a fixed ratio $r = \cF(\cB^\Omega) - D$ and then expressed in Theorem \ref{prop:CharDomOrder} the difference of the number of faces between a stuffed Walsh maps and a projected tree in terms of circuit ranks and genera. Yet we do not know if the ratio obtained for projected trees is the maximum. If true, it requires $\Omega$ to be an optimal pairing, so that $\cF(\cB^\Omega)$ is maximized.

Using the formula of Theorem \ref{prop:CharDomOrder} we have nevertheless proved that projected trees are dominant maps in some examples in Section \ref{sec:Applications} (projected trees are in fact dominant maps in all known cases), and identified the other dominant maps. This is an important step towards understanding and solving models with non-Gaussian behaviors as proposed in \cite{New1/N} such as that of Section \ref{sec:6VertexD=4}. The main criterion we have used is to minimize the number of \enquote{rainbow} cycles (whose edges do not all have the same color sets) which can be done using spanning forests.

Another interesting question is whether stuffed Walsh maps satisfy the topological recursion. Indeed, stuffed maps, introduced in \cite{Stuffed}, do satisfy the topological recursion. It is a recursion on the generating functions with $n$ marked faces and genus $g$, which comes out as a solution of Tutte's equations for combinatorial maps. In the case of stuffed maps, the recursion in fact receives a non-universal, stuffing-dependent part.

Finally, we would like to point out a crucial difference between the quartic case and the generic case presently studied. We have mentioned in Section \ref{sec:MatrixModel} that tensor integrals generate colored graphs and matrix integrals generate maps. This holds via a formal expansion of the integrals in powers of the counting parameters. However, the integrals themselves are expected to be analytic in a certain domain of the complex plane for the counting parameter (typically a cardioid) which in particular contains zero only on its boundary, meaning that there is more to the integrals than their formal expansions as power series in the counting parameters. In the quartic case, the bijection in fact lifts to a true equality between the tensor and matrix integrals, said to be a \emph{non-perturbative} equality. In our case however, the equality of Theorem \ref{thm:intmat1} only holds perturbatively, i.e. order by order in powers of the counting parameter. A non-perturbative version would be very important, but we cannot expect it to exist for arbitrary bubbles.

\section*{Acknowledgements}

The authors wish to thank A. Bres, S. Dartois and R. Gurau for useful discussions.

\appendix

\section{Edge-colored maps} \label{sec:EdgeColoredMaps}

\begin{prop}
Each pairing $\Omega$ of the vertices of $\cB\in\bB$ generates an infinite family of maps $\bM(\cB,\Omega)$ with boundary graph $\cB$. Since any map in $\partial^{-1}(\cB)$ induces a unique pairing, this defines a partition $\partial^{-1}(\cB)=\sqcup_{\Omega_{i}}\bM(\cB, \Omega_{i})$, where the union is over inequivalent pairings.
\end{prop}

\prf Any map has a unique boundary graph and induces a unique pairing of its vertices, so that the only unproven point is that given any $\cB$, the family $\bM_{\cB, \Omega}$ is infinite, which is true e.g. because the addition of degree 2 vertices on edges does not change the boundary graph.
\qed
 
Note that two pairings may generate the same family in which case they are equivalent. The more a bubble has symmetries, the fewer disjoint families there are. See for example the case of the complete bipartite graph in Subsection \ref{sec:K33}.
 
\begin{definition}
The \emph{power} $\delta(\cM)$ of an edge-colored map $\cM\in\bM_p$ is 
\begin{equation}
\delta(\cM) = \cF(\cM) - (D-1)\cE(\cM).
\end{equation}
\end{definition} 
It corresponds to the exponent of $N$ in the amplitude of $\cM$ in the Feynman expansion of the quartic melonic tensor model after a Hubbard-Stratonovich transformation \cite{DSDartois, BeyondPert}.

\begin{prop}
\label{prop:PowVacQuart}
The power of a map $\cM\in\bM_0$ is
\be
\label{eqref:PowerVacQuart}
\delta(\cM) = D(1-l(\cM))+2\sum_{i=1}^D \bigl(l(\cM^{(i)})-g(\cM^{(i)})\bigr),
\ee
where $l=\cE-\cV+1$ is the circuit rank of a connected graph, $g$ the genus of a map, given by $2-2g=\cV-\cE+\cF$, and $\cM^{(i)}$ is the monochromatic submap of color $i$.
\end{prop}
 
\prf By induction, it is easy to see that trees in $\bM_0$ have power $D$, hence the number of faces of a tree $\cT$ with $\cE(\cT)$ edges is
\begin{equation}
\cF(\cT) = (D-1)\cE(\cT) + D.
\end{equation}
One can then apply Theorem \ref{prop:CharDomOrder} to the quartic melonic case and notice that stuffed and projected maps are the same in this case, so that their circuit ranks are equal.
\qed

\begin{coroll}
\label{coroll:CorollMajTree}
The power $\delta$ of a map $\cM_0\in\bM_0$ verifies
\be
\label{eqref:CorollMajTree}
\delta(\cM_0)\le D,
\ee
and the inequality is saturated if and only if the map $\cM_0$ is a tree.
\end{coroll}
 
This is a very well-known result for colored graphs \cite{Uncolored}, see also \cite{DSDartois, BeyondPert} for his translation in terms of edge-colored maps. In the context of edge-colored maps, there is a particularly simple proof: by induction, it is shown that $\delta(\cM)=D$ for trees, and since every edge carries a single color, $|\cI_2(e)| \leq 1$ so that Lemma \ref{lemma:Unhooking} ensures that the number of faces increases upon unhooking edges until $\cM$ becomes a tree. We offer here an alternative proof using circuit ranks.

\prf To prove that result, we use \eqref{eqref:PowerVacQuart} and first show that
\be \label{ColoredCycleBound}
\frac{D}{2}l(\cM_0)\ge\sum_{i=1}^Dl(\cM_0^{(i)}).
\ee
For each $i\in[D]$, consider a forest $\cT^{(i)}$ spanning $\cM_0^{(i)}$. There are $l(\cM_0^{(i)})$ edges in $\cM^{(i)}_0\setminus \cT^{(i)}$ and they identify fundamental cycles. The union $\bigcup_{i\in[D]} \cT^{(i)}$ is connected and spans $\cM_0$. A spanning tree $\cT$ of $\bigcup_{i\in[D]} \cT^{(i)}$ is thus a spanning tree of $\cM_0$. Since all the edges in $\bigcup_{i\in[D]} (\cM_0^{(i)}\setminus \cT^{(i)})$ are distinct, 
\begin{equation}
l(\cM_0)\ge\sum_{i=1}^Dl(\cM_0^{(i)}),
\end{equation}
and \eqref{ColoredCycleBound} follows from $D\ge2$. Since $g(\cM_0^{(i)})\ge0$, one obtains \eqref{eqref:CorollMajTree}. 

Furthermore, if $e\not\in\cT$, then either $e\in \cup_{i\in[D]} \cT^{(i)}$ or $e\not\in\cup_{i\in[D]} \cT^{(i)}$ and since $\cT$ is a spanning tree for $\cup_{i\in[D]}\cT^{(i)}$, we find
\be
l(\cM_0)=l(\cup_{i=1}^D\cT^{(i)})+\sum_{i=1}^Dl(\cM_0^{(i)}),
\ee
so that \eqref{eqref:PowerVacQuart} can be rewritten as 
\be
\delta(\cM_0)=D-Dl(\cup_{i=1}^D\cT^{(i)})-(D-2)\sum_{i=1}^Dl(\cM_0^{(i)})-2\sum_{i=1}^Dg(\cM_0^{(i)}),
\ee
and the inequality is saturated if and only if $l(\cup_{i=1}^D\cT^{(i)})=l(\cM_0^{(i)})=g(\cM_0^{(i)})=0$, i.e. if and only if $\cM_0$ is a tree.
\qed
 
\begin{lemma}
\label{lemma:VacVersusExt}
Given an edge-colored map $\cM_p$ with $p$ cilia and $\cM_0$ obtained from $\cM_p$ by removing the cilia and merging their incident corners, the power of $\cM_p$ is related to that of $\cM_0$ by the following relation,
\be
\delta(\cM_p)=\delta(\cM_0)-\cF(\cB^\Omega)
\ee
where $\cB = \partial \cM_p$, $\Omega$ is the pairing of $\cB$ induced by $\cM_p$ and $\cB^\Omega$ is the corresponding covering.
\end{lemma}

\prf $\delta(\cM_p)$ and $\delta(\cM_0)$ only differ due to the numbers of faces. Every face of $\cM_p$ is a face of $\cM_0$. However, removing the cilia and merging their incident corners turn the broken faces into a set of closed faces. This operation consists in adding an edge of color $0$ between the white vertex and the black vertex of each pair defined by $\Omega$. Those faces are thus described by the covering $\cB^{\Omega}$.
\qed
 
\begin{lemma}
 \label{lemma:IneqRain}
Let $\cM$ be an edge-colored map with boundary graph $\partial \cM = \cB$ and induced pairing $\Omega$ and monochromatic submaps denoted $\cM^{(i)}$.
The directed graph $\cB_{\circlearrowleft,\Omega}$ is obtained from $\cB$ by orienting the edges from white to black vertices, then contracting the pairs of vertices defined by $\Omega$ and removing the edges which connect the vertices of a pair. Then 
 \be
 \label{eqref:IneqOpt}
 l(\cM)-\sum_{i=1}^D l(\cM^{(i)})\ge l(\cB_{\circlearrowleft,\Omega})-\sum_{i=1}^D l(\cB_{\circlearrowleft,\Omega}^{(i)})=\opt_{\cB}(\Omega).
 \ee
We call $\opt_{\cB}(\Omega)\in\bN$ the optimality of the pairing $\Omega$.
 \end{lemma}
 
The quantity $l(\cM)-\sum_{i=1}^D l(\cM^{(i)})$ is the number of independent (proper) \emph{rainbow} cycles of a graph, i.e. cycles in which edges do not all have the same color. When each monochromatic submap is a forest, it coincides with the circuit rank of the map itself.
 
\prf For each color $i$, consider $\cT^{(i)}$ a forest that spans $\cM^{(i)}$. As in the proof of Corollary $ \ref{coroll:CorollMajTree}$, the number of rainbow cycles of $\cM$ is
\be
l(\cM)-\sum_{i=1}^Dl(\cM^{(i)})=l(\cup_{i=1}^D\cT^{(i)}).
\ee
Consider $\cT_{ext}^{(i)}$ the subforest of $\cT^{(i)}$ made of its connected components which have at least {\it two} marked vertices. The union $\cup_{i=1}^D\cT_{ext}^{(i)}$ is a subgraph of $\cup_{i=1}^D\cT^{(i)}$, hence
\be
l(\cup_{i=1}^D\cT^{(i)})\ge l(\cup_{i=1}^D\cT_{ext}^{(i)}).
\ee
Notice that every cilium belongs to at least one $\cT_{ext}^{(i)}$, else there would be one cilium for which all broken faces which start at this cilium end at it and $\cB$ would be disconnected. Therefore $\cup_{i=1}^D\cT_{ext}^{(i)}$ is connected. Moreover, $\cE(\cup_{i=1}^D\cT_{ext}^{(i)})=\sum_{i=1}^D\cE(\cT_{ext}^{(i)})$, so that
\be
l(\cup_{i=1}^D\cT_{ext}^{(i)})=\sum_{i=1}^D\cE(\cT_{ext}^{(i)})-\cV(\cup_{i=1}^D\cT_{ext}^{(i)})+1.
\ee
Since $\cT_{ext}^{(i)}$ is a forest, $\cE(\cT_{ext}^{(i)})=\cV(\cT_{ext}^{(i)})-k(\cT_{ext}^{(i)})$. Furthermore, a vertex $v\in\cup_{i=1}^D \cT_{ext}^{(i)}$ is in as many distinct subforests $\cT_{ext}^{(i_1)}, \cT_{ext}^{(i_2)}, \dotsc$ as the number $\col(v)$ of incident colors, meaning
\begin{equation}
\sum_{i\in[D]} \cV(\cT_{ext}^{(i)}) = \sum_{v\in \cup_{i=1}^D\cT_{ext}^{(i)}} \bigl(\col(v)-1\bigr).
\end{equation}
It comes that
\be
\label{eqref:CyclExt}
l(\cup_{i=1}^D\cT_{ext}^{(i)})
=\sum_{v\in\cup_{i=1}^D\cT_{ext}^{(i)}}(\col(v)-1)-\sum_{i=1}^Dk(\cT_{ext}^{(i)})+1.
\ee
Moreover, $\cup_i\cT_{ext}^{(i)}$ has no isolated vertices, hence $\col(v)-1\ge0$, and
\be
\sum_{v\in\cup_{i=1}^D\cT_{ext}^{(i)}}(\col(v)-1)\ge\sum_{\substack{ \text{$v$ ciliated}\\\text{in $\cup_{i=1}^D\cT_{ext}^{(i)}$}}}(\col(v)-1)=\sum_{v\in\cB_{\circlearrowleft,\Omega}}(\col(v)-1).
\ee
To understand the last equality, construct the graph $\cG_{\circlearrowleft, \cup_{i}\cT_{ext}^{(i)}}$ as in Section \ref{sec:BoundaryGraphs}: only keep the ciliated vertices and draw a directed edge for each broken face. Since all connected components of $\cT_{ext}^{(i)}$ have at least two ciliated vertices, $\cG_{\circlearrowleft, \cup_{i}\cT_{ext}^{(i)}}$ has no loop-edges. That implies $\cG_{\circlearrowleft, \cup_{i}\cT_{ext}^{(i)}} = \cB_{\circlearrowleft,\Omega}$. In addition, since the closed walks of fixed color $j$ and which encounter cilia on $\cup_{i=1}^D\cT_{ext}^{(i)}$ are in one-to-one correspondence with the connected components of $\cB_{\circlearrowleft,\Omega}^{(j)}$, and since there is at least one such walk per connected component of $\cT_{ext}^{(j)}$, it comes that $k(\cB_{\circlearrowleft,\Omega}^{(i)})\ge k(\cT_{ext}^{(i)})$. We therefore get the bound
\be
\label{eqref:PrfIneq}
l(\cup_{i=1}^D\cT_{ext}^{(i)})\ge\sum_{v\in\cB_{\circlearrowleft,\Omega}}(\col(v)-1)-\sum_{i=1}^Dk(\cB_{\circlearrowleft,\Omega}^{(i)})+1.
\ee
To conclude, we use the fact that $\cB_{\circlearrowleft,\Omega}^{(i)}$ consists of disjoint cycles, so that: $\sum_{v\in\cB_{\circlearrowleft,\Omega}}(\col(v)-1)=\cE(\cB_{\circlearrowleft,\Omega})-\cV(\cB_{\circlearrowleft,\Omega})$, and $l(\cB_{\circlearrowleft,\Omega}^{(i)})=k(\cB_{\circlearrowleft,\Omega}^{(i)})$.
\qed

\begin{lemma}[Optimality]\label{lemma:OptPow}
The number of faces of the covering $\cB^\Omega$ of $\cB$ a bubble with $\cV$ white vertices is related to the optimality of the pairing as follows,
\be
\cF(\cB^\Omega) = 1 + (D-1)\cV - \opt_{\cB}(\Omega).
\ee 
\end{lemma}

Recall from Definition \ref{def:OptimalPairing} that an optimal pairing is one which maximizes $\cF(\cB^\Omega)$. It is thus equivalent to minimizing the optimality.

\prf This relies on the fact that the faces of color $i$ of $\cB^\Omega$ of length two (an edge of color $i$ and an edge of color 0) are the isolated vertices of $\cB_{\circlearrowleft,\Omega}^{(i)}$ while the the faces of length greater than two are in one-to-one correspondence with the cycles of $\cB_{\circlearrowleft,\Omega}^{(i)}$. Each vertex of $\cB_{\circlearrowleft,\Omega}$ is an isolated vertex in $(D - \col(v))$ monochromatic submaps. Thus,
\be
\cF(\cB^\Omega) = \sum_{i=1}^D \cF_{i}(\cB^\Omega) = \sum_{v\in\cB_{\circlearrowleft,\Omega}}(D-\col(v)) + \sum_{i=1}^D l(\cB_{\circlearrowleft,\Omega}^{(i)}) = D\cV - \cE(\cB_{\circlearrowleft,\Omega}) + \sum_{i=1}^D l(\cB_{\circlearrowleft,\Omega}^{(i)}).
\ee 
Therefore
\be
1 - \cF(\cB^\Omega) + (D-1)\cV = 1 - \cV + \cE(\cB_{\circlearrowleft,\Omega}) - \sum_{i=1}^D l(\cB_{\circlearrowleft,\Omega}^{(i)}) = l(\cB_{\circlearrowleft,\Omega}) - \sum_{i=1}^D l(\cB_{\circlearrowleft,\Omega}^{(i)}) = \opt_{\cB}(\Omega),
\ee
using the definition of the optimality \eqref{eqref:IneqOpt}.
\qed

\begin{theorem}
Let $\cM_p\in \bM_p$ be an edge-colored map, with boundary graph $\cB = \partial \cM_p$, $p$ cilia (number of white vertices of $\cB$). Then its power is bounded as follows,
\be
\label{PowDomB}
\delta(\cM_p)\leq -(D-1)(p+\opt_{\cB}(\Omega)-1).
\ee
Moreover the equality holds if and only if 
\begin{itemize}
\item $\cM_p^{(i)}$ is a forest for each $i\in[D]$, and
\item denoting $\cT_{ext}^{(i)}$ the restriction of $\cM_p^{(i)}$ to its trees which have at least two cilia, those trees can meet on ciliated vertices only, and
\item $\cM_p\setminus \bigcup_{i\in[D]} \cT_{ext}^{(i)}$ is a forest whose trees each meet $\bigcup_{i\in[D]} \cT_{ext}^{(i)}$ at a single vertex.
\end{itemize}
\end{theorem}

\prf We start with Lemma \ref{lemma:VacVersusExt}, $\delta(\cM_p) = \delta(\cM_0) - \cF(\cB^\Omega)$ and apply Proposition \ref{prop:PowVacQuart} to $\delta(\cM_0)$ and Lemma \ref{lemma:OptPow} to $\cF(\cB^\Omega)$. That gives
\begin{equation}
\delta(\cM_p) = (D-1)(1-p) -D l(\cM_0) + 2 \sum_{i\in [D]} l(\cM_0^{(i)}) - g(\cM_0^{(i)}) + \opt_\cB(\Omega).
\end{equation}
Rewrite
\begin{equation}
-D l(\cM_0) + 2 \sum_{i\in [D]} l(\cM_0^{(i)}) = -D\Bigl( l(\cM_0) - \sum_{i\in [D]} l(\cM_0^{(i)})\Bigr) - (D-2) \sum_{i\in [D]} l(\cM_0^{(i)}).
\end{equation}
Since $(D-2) l(\cM_0^{(i)}) \geq0$ as well as $g(\cM_0^{(i)}) \geq0$ we get the bound
\begin{equation}
-D l(\cM_0) + 2 \sum_{i\in [D]} l(\cM_0^{(i)}) - g(\cM_0^{(i)}) \leq -D\Bigl( l(\cM_0) - \sum_{i\in [D]} l(\cM_0^{(i)})\Bigr) \leq -D \opt_\cB(\Omega),
\end{equation}
where the last inequality is due to Lemma \ref{lemma:IneqRain}. This proves \eqref{PowDomB}.

One gets the equality in \eqref{PowDomB} if and only if $\cM_p$ verifies for all $i\in[D]$
\be
\label{eqref:Cond1}
l(\cM_p^{(i)})=g(\cM_p^{(i)})=0,
\ee
i.e. $\cM_p^{(i)}$ is a forest, and 
\be
\label{eqref:Cond2}
l(\cM_p)-\sum_{i=1}^Dl(\cM_p^{(i)})=\opt_{\cB}(\Omega).
\ee
Using the notations and results of the proof of Lemma \ref{lemma:IneqRain}, we have the sequence of bounds
\begin{equation}
\label{eqref:SummIneqRain}
\begin{aligned}
l(\cM_p)-\sum_{i=1}^Dl(\cM_p^{(i)}) = l(\cup_i\cT^{(i)}) \geq l(\cup_i\cT_{ext}^{(i)}) &= \sum_{v\in\cup_i\cT_{ext}^{(i)}}(\col(v)-1)-\sum_{i=1}^D k(\cT_{ext}^{(i)})+1 \\
&\geq \sum_{\substack { \text{$v$ ciliated}\\\text{in $\cup_i\cT_{ext}^{(i)}$}}}(\col(v)-1)-\sum_{i=1}^Dk(\cB_{\circlearrowleft,\Omega}^{(i)})+1 = \opt_\cB(\Omega),
\end{aligned}
\end{equation}
so that \eqref{eqref:Cond2} is satisfied if and only if every inequality in \eqref{eqref:SummIneqRain} is an equality. This is equivalent to
\begin{itemize}
\item the connected components of $\cT_{ext}^{(i)}$ are in one-to-one correspondence with those of $\cB_{\circlearrowleft,\Omega}^{(i)}$,
\item the unciliated vertices of $\cup_i\cT_{ext}^{(i)}$ are reached by a single color,
\item $\cup_i\cT^{(i)}\setminus\cup_i\cT_{ext}^{(i)}$ is a forest whose connected components each share a single vertex with $\cup_i\cT_{ext}^{(i)}$.
\end{itemize}
Equation \eqref{eqref:Cond1} makes the first condition trivial, and reduces the other two to the statements listed in the theorem.
\qed

This result leads to the notion of \emph{optimal maps}, which are maps satisfying $\delta(\cM_p)= -(D-1)(p+\opt_{\cB}(\Omega)-1)$ where $\Omega$ is an optimal pairing of $\cB$, i.e. $\opt_\cB(\Omega)=\min_{\Omega_\alpha} \opt_{\cB}(\Omega_\alpha)$.

\section{Recovering the Hubbard-Stratonovich transformation for the quartic melonic model}
\label{HubStrat}

We want to show the equality between the matrix model of Theorem \ref{thm:intmat1} and the matrix model \cite{BeyondPert} obtained by the Hubbard-Stratonovich transformation in the case of quartic melonic bubbles. Restricting to the melonic bubble $\cB_i$, Theorem \ref{thm:intmat1} gives
\be
\label{eqref:ZQuartMel}
 Z_{\cB_i}(\lambda_i, N)
=\int e^{ -  \Tr\bigl[ {N^{D-1}}\lambda_i\sigma_i^2-\ln(\un^{\otimes D}+{\bf \bar \fsig _i})  \bigr]}d\mu_{0}(\sigma_i,{\bf \bar  \sigma_i}).
\ee
Its perturbative expansion can be performed by seeing the Gaussian measure as a differential operator,
\begin{align}
Z_{\cB_i}(\lambda_i, N)
&=\biggl[ \exp\biggl\{ \frac{1}{N^{D-1}}\sum_{a,b=1}^N\frac{\partial}{\partial\sigma_{i\mid a,b}}\frac{\partial}{\partial\bar\sigma_{i\mid a,b}}  \biggr\}     e^{ -  \Tr\bigl[ {N^{D-1}} \lambda_i\sigma_i^2-\ln(\un^{\otimes D}+{\bf \bar \fsig _i})  \bigr]}\biggr]_{\sigma_i=0}\\
&
\label{eqref:ExprZ}
=\sum_{p,q\ge0} \frac{(-\lambda_i)^qN^{(D-1)(q-p)}}{p!q!}
\biggl[\biggl( \sum_{a,b=1}^N\frac{\partial}{\partial\sigma_{i\mid a,b}}\frac{\partial}{\partial\bar\sigma_{i\mid a,b}}  \biggr)^p     \biggl( \sum_{c,d=1}^N\sigma_{i\mid c,d}\bar\sigma_{i\mid c,d}\biggr)^q e^{-\Tr\ln(\un^{\otimes D}+{\bf \bar \fsig _i})  }\biggr]_{\sigma_i=0}.
\end{align}

Denoting $F(\bar\fsig_i)=e^{-\Tr\ln(\un^{\otimes D}+{\bf \bar \fsig _i})}$ and $\binom{p}{\{p_{ab}\}}=\frac{p!}{\prod_{a,b}p_{ab}!}$, the expression under brackets rewrites
\bea
&S_{pq}
&=\sum_{  \substack{ {\{p_{ab}\}\mid\sum p_{ab}=p }\\{ \{q_{ab}\}\mid\sum q_{ab}=q } }  } \binom{p}{\{p_{ab}\}}\binom{q}{\{q_{ab}\}}\biggl[ \prod_{a,b=1}^N \biggl(\frac{\partial}{\partial\sigma_{i\mid a,b}}\frac{\partial}{\partial\bar\sigma_{i\mid a,b}} \biggr)^{p_{ab}}  \bigl(\sigma_{i\mid a,b}\sigma_{i\mid b,a}\bigr)^{q_{ab}} F(\bar\fsig_i) \biggr]_{\sigma_i=0}.\nonumber
\eea
Now noticing that for any set of quantities $\{A_{ab}\}$, 
\be
\label{eqref:ProdSym}
\prod_{a,b}(A_{ab}A_{ba})^{q_{ab}}=\prod_{a,b}(A_{ab})^{q_{ab}+q_{ba}},
\ee
this becomes, replacing $A_{ab}$ with $\sigma_{i\mid a,b}$ and evaluating the derivatives over $\sigma_{i\mid a,b}$,
\begin{equation}
\begin{aligned}
S_{pq} &= \sum_{  \substack{ {\{p_{ab}\}\mid\sum p_{ab}=p }\\{ \{q_{ab}\}\mid\sum q_{ab}=q } }  } \binom{p}{\{p_{ab}\}}\binom{q}{\{q_{ab}\}}
\biggl[ \prod_{a,b=1}^N p_{ab}!\ \delta_{p_{ab}}^{q_{ab}+q_{ba}} \biggl(\frac{\partial}{\partial\bar\sigma_{i\mid a,b}}\biggr)^{p_{ab}}F(\bar\fsig_i) \biggr]_{\sigma_i=0} \\
&=\sum_{ \{q_{ab}\}\mid\sum q_{ab}=q} \binom{q}{\{q_{ab}\}} \biggl[ \prod_{a,b=1}^N \ \biggl(\frac{\partial}{\partial\bar\sigma_{i\mid a,b}}\frac{\partial}{\partial\bar\sigma_{i\mid b,a}}\biggr)^{q_{ab}}F(\bar\fsig_i) \biggr]_{\bar\sigma_i=0}p!\ \delta_{p}^{2q} \\
&=\biggl[ \biggl(\sum_{a,b=1}^N   \frac{\partial}{\partial\bar\sigma_{i\mid a,b}}\frac{\partial}{\partial\bar\sigma_{i\mid b,a}} \biggr)^q F(\bar\fsig_i) \biggr]_{\bar\sigma_i=0}p!\ \delta_{p}^{2q}
\end{aligned}
\end{equation}
where we used that $\sum_{ \{p_{ab}\}\mid\sum p_{ab}=p}p!\prod_{a,b=1}^N\delta_{p_{ab}}^{q_{ab}+q_{ba}}=p!\ \delta_{p}^{2q}$, and again relation \eqref{eqref:ProdSym} applied to $\frac{\partial}{\partial\bar\sigma_{i\mid a,b}}$.

We have thus completely gotten rid of the matrix $\sigma_i$ and are left with $\bar{\sigma}_i$. Inserting this back into \eqref{eqref:ExprZ}, the partition function finally writes
\begin{equation}
\begin{aligned}
Z_{\cB_i}(\lambda_i, N)&=\biggl[ \exp\biggl\{-\frac{\lambda_i}{ N^{(D-1)}}\sum_{a,b=1}^N \frac{\partial}{\partial\bar\sigma_{i\mid a,b}}\frac{\partial}{\partial\bar\sigma_{i\mid b,a}}   \biggr\} F(\bar\fsig_i) \biggr]_{\bar\sigma_i=0}\\
&=\int e^{ +\frac{\lambda_i}{N^{D-1}}\Tr\bar \sigma_i^2}F(\bar\fsig_i)\frac{d{\bar  \sigma_i}}{\sqrt{\pi}^{N^2}}
\end{aligned}
\end{equation}
and therefore
\begin{equation}
Z_{\cB_i}(\lambda_i, N) = i\biggl(\frac{N^{D-1}}{2\lambda_i}\biggr)^{\frac{1}{2}}\int e^{ -\frac{1}{2}\Tr \sigma_i^2}F\biggl(i\biggl(\frac{N^{D-1}}{2\lambda_i}\biggr)^{\frac{1}{2}}\fsig_i\biggr)\frac{d{\sigma_i}}{\sqrt{\pi}^{N^2}},
\end{equation}
as expected from the Hubbard-Stratonovich transformation \cite{BeyondPert}.


\end{document}